\documentclass[a4paper,reqno,11pt]{amsart}
\usepackage{amsfonts,latexsym,amsmath,showkeys}
\newtheorem{theorem}{Theorem}[section]
\newtheorem{lemma}[theorem]{Lemma}
\newtheorem{proposition}[theorem]{Proposition}

\newtheorem{hypothesis}[theorem]{Hypothesis}

\theoremstyle{definition}
\newtheorem{definition}[theorem]{Definition}
\newtheorem{example}[theorem]{Example}

\theoremstyle{remark}
\newtheorem{remark}[theorem]{Remark}

\numberwithin{equation}{section}

\def\Dim{\noindent\hbox{{\bf Proof.}$\;\; $}}          
\def\finedim{{\hfill\hbox{\enspace${ \square}$}} \smallskip}    
\def\sqr#1#2{{\vcenter{\vbox{\hrule height .#2pt
     \hbox{\vrule width .#2pt height#1pt \kern#1pt \vrule
     width .#2pt} \hrule height .#2pt}}}}
\def\square{\mathchoice\sqr54\sqr54\sqr{4.1}3\sqr{3.5}3}

\def\media{\mathbb E}           
              \def\E{\mathbb E}       
             \def\R{\reali}
             \def\reali{\mathbb R}   

\def\F{{\mathcal F}}
 \def\B{{\mathcal B}}          \def\K{\mathcal K}

\def\P{{\mathcal P}}

\def\Et{\mathbb{E}^{\F_t}}

\def\<{\langle}
\def\>{\rangle}
\begin{document}
$ \ $ \\
\bigskip

\centerline{\Large }
\centerline{\Large }
\bigskip

\title[Mild Solution of BSREs]{Well Posedness of Operator Valued Backward Stochastic Riccati  Equations in Infinite Dimensional Spaces}
\medskip

\maketitle
\centerline{Giuseppina Guatteri}

\centerline{ Dipartimento di Matematica,
  Politecnico di Milano}
\centerline{piazza Leonardo da Vinci 32, 20133 Milano, Italia}
\centerline{e-mail: guatteri@mate.polimi.it}
\bigskip

\centerline{Gianmario Tessitore}
\centerline{ Dipartimento di Matematica e Applicazioni,
Universit\`{a} di Milano-Bicocca,}
\centerline{via R. Cozzi 53 - Edificio U5, 20125 Milano, Italia}
\centerline{e-mail: gianmario.tessitore@unimib.it}
\date{\today}

\email{}

\date{\today}

\bibliographystyle{plain}

\begin{abstract}
We prove existence and uniqueness of the mild solution of an infinite dimensional, operator valued, backward stochastic Riccati equation. We exploit
the regularizing properties of the semigroup generated by the unbounded operator involved in the equation. Then the results will be applied to characterize the value function and optimal feedback law for a infinite dimensional, linear quadratic control problem with stochastic coefficients.
\end{abstract}

\smallskip {\bf Key words.} 
Backward Stochastic Differential Equations in infinite dimensions, Riccati equation, linear quadratic optimal control, Hilbert spaces, stochastic coefficients.

 \smallskip {\bf AMS.} 
93E20, 60H10, 49A60,35R60


\section{Introduction}

The present paper is concerned with the following infinite dimensional Backward Stochastic Riccati Equation
(BSRE)
\begin{equation}\label{Riccati-intro}
\left\{  \begin{array}{l}   -dP_t= (A' P_t + P_t A +C_t'Q_t + Q_t C_t  + C'_t P_t C_t- P_t B_tB'_tP_t + S_t  ) \,dt   - Q_t d\, W_t, \\ \\
\quad P_T=M
  \end{array}
\right .
\end{equation}
where $A$ is a self adjoint operator on the Hilbert space $H$  generating the analytic semigroup $(e^{tA})$;  $(W_t)_{t\geq 0}$ is a real valued standard Brownian motion;  $(B_t)$, $(C_t)$, $(S_t)$ are
operator valued adapted processes.
 The unknown of the equation is the couple $(P,Q)$ of operator valued adapted processes .

As it is well known see \cite{Yong_Zhou}  the   above equation represents the value function of
a linear quadratic optimal control problem
 involving a Hilbert valued state equation with stochastic coefficients (in particular of a control
 problem with evolution modelled  by a parabolic SPDE with stochastic coefficients).
It is also well known that, as soon as the solution of the BSRE is obtained,
then the synthesis of the optimal control easily follows with a clear applicative interest.

Moreover the special case in which $B_t\equiv 0 $ (the so called Lyapunov equation)
turns out to be  essential  in the formulation of the Pontyagin maximum principle
for controlled systems described by stochastic partial differential equations
(see \cite{LiTang} \cite{LuZhang}, \cite{Du-Meng}, \cite{Fu-Hu-Te-1}, \cite{Fu-Hu-Te-2}).
This in particular happens in the so called general case in which the space of controls is not convex and the control affects the diffusion term as well (see \cite{Pe99}). Indeed this is the case in which the second variation process, that satisfies an operator Lyapunov equation, has to be introduced. In this context the research on backward evolution equations in spaces of linear operators has gained recently a relevant interest.

The study of BSREs in finite dimensional spaces had quite a long story between the pioneering paper
by J.M. Bismut and then S. Peng (see \cite{Bi76} and \cite{Peng}) and the conclusive paper by S. Tang (see
\cite{Tang}) where existence and uniqueness is proved in the most general case.

On the contrary the study of BSREs in infinite dimensional spaces adds specific new difficulties and few results are available.
As far as the Lyapunov equation is concerned in \cite{LiTang} the solution is obtained when the final condition $M$ and the forcing term $S$ are Hilbert-Schmidt operators
(condition that is rarely satisfied) while in \cite{Du-Meng}, \cite{Fu-Hu-Te-2} the process $P$ is  characterized by an energy equality
involving a suitable forward stochastic differential equation in $H$. Finally in \cite{LuZhang} the concept of \emph{transposed} solution is given which again consists in a characterization of $P$ and $Q$ by a  suitable   \emph{duality} relation that involves an infinite dimensional  forward equation. We notice that in all the above cases no explicit differential or integral equation directly satisfied by $P$ and $Q$ is presented.

Regarding the Riccati equation (that, differently from the Lyapunov equation, is non linear) in \cite{GuaTess}
we proposed to characterize the $P$-part of the solution using the concept of strong solution which is of common  use in PDE theory (see \cite{BeDPDeMi} or \cite{Lunardi}). Roughly speaking we characterize the solution as the limit of a sequence of equations with regular (in this case Hilbert-Schmidt) data. This result is good enough to be applied to the corresponding linear quadratic control problem  but has the drawback of not saying anything on the martingale term of the solution (the $Q$-term) and consequently not giving the representation through a (differential) equation.

The origin of the difficulties to deal with stochastic backward Riccati (or even Lyapunov) equation in the infinite dimensional case is in the fact that the natural space in which it should be treated is the space $L(H)$ of bounded linear operators in $H$ which is only a Banach space that does not enjoy any of the regularity properties (as UMD or M-type condition) allowing to establish an analogue of the classical Hilbertian stochastic calculus. Moreover although, as we have said above, different characterization of the solution have been recently proposed, it seems to us that the natural notion of solution is the one of \textit{mild} solution introduced in the theory of infinite dimensional BSDEs since the seminal paper by \cite{HuPeng1991}. We finally notice that this way both the $P$ and the $Q$ part of the unknown is characterized by a differential equation.

As far as we know this is the first paper in which existence and uniqueness of a mild solution of equation (\ref{Riccati-intro}) is obtained.  Indeed we show that $(P,Q)$ is the unique couple of processes (with suitable regularity) verifying 
\begin{align}
P(t)&=e^{(T-t)A'}Me^{(T-t)A}+ \int_t^Te^{(s-t)A'}
S(s)e^{(s-t)A}\,ds\nonumber \\
&+ \int_t^Te^{(s-t)A'}
\Big[C'(s)P(s)C(s)+C'(s)Q(s)+Q(s)C(s)\Big]e^{(s-t)A}\,ds
\\\nonumber
 &+ \int_t^Te^{(s-t)A^*} Q(s)e^{(s-t)A}\,dW(s)\quad\quad \mathbb{P}-\text{a.s.}
\end{align}
where $P$ is a predictable process with values in the space of bounded non negative, simmetric, linear operators in $H$ which as we said is, in some sense, the natural space for the equation. On the contrary the identification of the right operators space for the evolution of $Q$ is the main achievement of this work. We shall prove existence and uniqueness of $Q$ as a square-integrable, adapted, process in a space $\K $ of Hilbert-Schimidt operators from  suitable domains of the fractional powers of $A$ (see (\ref{defK})). This is an  Hilbert space, large enough to contain all bounded operators. This choice will allow to recover stochastic calculus tools. The price to pay is that the term
$C'_tQ_t + Q_t C  + C'_t P_t C_t
$ becomes unbounded on $\K$ . This difficulty will be handled exploiting in a careful (and non completely standard) way the regularizing properties of the semigroup generated by $A$. By the way we have to say that our results rely on the specific properties of $A$ that we assume to be self-adjoint with rapidly increasing eigenvalues. Nevertheless our assumptions can cover important classes of strongly elliptic differential operators.

The structure of the proof will be the following: first we introduce suitable approximations  of the equation
 (see (\ref{Lyapmildappr}) that can be treated bt the standard Hilbert-Schmidt theory. Then showing the needed convergence estimates we  prove existence and uniqueness of the solution to a simplified Lyapunov equation (see \ref{Gamma}). An a-priori estimate (see (\ref{apriori}) helps to prove convergence and gives uniqueness. Consequently a fixed point argument yields existence and uniqueness
of a solution to the Lyapunov equation. Finally, in Section \ref{Sec-LQ}, we exploit the interplay between the Riccati equation and the corresponding optimal control problem to obtain existence and uniqueness of the mild solution to the BSRE and the synthesis of the optimal control.

We notice that the optimal control problem is given by the
 following {\em state equation}:
\begin{equation}\label{stato.intro}
\left\{
\begin{array}{ll}
dy(t)=(Ay(t)+B(t) u(t)) \,dt + C(t) y(t) \,dW(t) & t \in [0,T]  \\
y(0)=x
\end{array}
\right.
\end{equation}
where $y$ is the \emph{state} of the system and
$u$ is the \emph{control}; $y$ and $u$ are adapted processes with values in  $H$,
and by the following
 quadratic \emph{cost functional}:
\begin{equation}\label{costo.intro}
\E \int_0^T \left(|\sqrt{S}_sy_s|^2 +
|u(s)|^2\right)\,ds + \E \langle My_T,y_T\rangle.
\end{equation}

\section{Main Notation and Assumptions}

%
%
%
%
%
%

\noindent{\bf Some classes of stochastic processes} \\
Let $G$ be any separable Hilbert space. By $\mathcal{P}$ we denote
the predictable $\sigma$-field on $ \Omega \times [0,T]$ and by
$\mathcal{B}(G)$ Borel $\sigma$-field on $G$.  The following
classes of processes will be used in this work
\begin{itemize}
\item $L^p_\P(\Omega \times [0,T];G)$, $p\in [1,+\infty]$ denotes
subset of  $L^p(\Omega \times [0,T];G)$, given by all equivalence
classes
 admitting a predictable version. This space is endowed with the natural
 norm.
\item $C_\P([0,T];L^p(\Omega;G))$ denotes the space of $G$-valued processes $Y$ such that
 $Y : [0,T] \to L^p(\Omega,G)$ is continuous
 and $Y$ has a predictable modification, endowed with the norm:
\begin{equation*}
|Y|^p_{C_\P([0,T];L^p(\Omega;G))}=\sup_{t \in [0,T]}\media
|Y_t|^p_G
\end{equation*}
Elements of $C_\P([0,T];L^p(\Omega;G))$ are identified up to
modification.
\item $L^p_\P(\Omega;C([0,T];G))$ denotes the space of predictable processes $Y$
with continuous paths in $G$, such that the norm
\begin{equation*}
|Y|^p_{L^p_\P(\Omega;C([0,T];G))}=\media \sup_{t \in
[0,T]}|Y_t|_G^p
\end{equation*}
is finite. Elements of this space are defined up to
indistiguishibility.
\end{itemize}

Now let us consider the space $L(G)$ of linear and bounded
operators from $G$ to $G$. This space,
as long as $G$ is infinite dimensional, is not separable, see
\cite[pag.23]{DPZ1}, therefore we introduce the following
$\sigma$-field:
\begin{equation*}
\mathcal {L}_S = \{ T \in L(G): Tu \in A  \}, \text { where } u
\in G \text{ and } A \in \mathcal{B}(G)
\end{equation*}
Following again \cite{DPZ1} the elements of $\mathcal{L}_S$ are
called {\em strongly measurable}.

\medskip
We notice that the maps $P \to |P|_{L(G)}$ and $(P,u) \to Pu$ are
measurable from $(L(G),\mathcal {L}_S)$ to $\mathbb{R}$ and from
$(L(G)\times G, \mathcal{L}_S\otimes\mathcal{B}(G))$ to
$(G,\mathcal{B}(G)) $ respectively.
\medskip

Moreover $\mathcal {L}_S$ is equivalent to the weak
$\sigma$-field:
\begin{equation*}
\mathcal {L}_S = \{ T \in L(G): (Tu,x) \in A  \}, \text { where }
u,x \in G \text{ and } A \in \mathcal{B}(\mathbb{R})
\end{equation*}

 We define the following spaces:
\begin{itemize}
\item $L^{\infty}_{\P,S}((0,T)\times \Omega; L(G))$
a space of predictable processes $Y$ from $(0,T)$ to $L(G)$,
endowed with the $\sigma$-field $\mathcal{L}_S$. For each element
$Y$ there exists a constant $C
>0$ such that:
\begin{equation*}
 |Y(t,\omega)|_{L(G)}  \leq C \quad\quad\quad \mathbb{P}-\text{a.s. for a.e. } t \in (0,T)
\end{equation*}
 In the same way we define
$L^\infty_S(\Omega,\F_T;L(G))$ as the space of maps $Y$ from
$(\Omega,\F_T)$ into $(L(G),\mathcal{L}_S)$ such that there exists
a positive constant $K$ such that:
\begin{equation*}
|Y(\omega)|_{L(G)} \leq K \qquad\quad \mathbb{P}-\text{a.s.}
\end{equation*}
\end{itemize}
Elements of this space are identified up to modification.

By $\Sigma(G)$ we denote the subspace of  all symmetric and operators and by $\Sigma^+(G)$  the convex subset of all
positive semidefinite operators.
 We
define identically the following spaces:
$L^{\infty}_{\P,S}((0,T)\times \Omega; \Sigma^+(G))$,
$L^1_{\P,S}((0,T);L^\infty(\Omega,\Sigma^+(G)))$
 and $L^\infty_S(\Omega,\F_T;\Sigma^+(G))$.

$ $

\noindent{\bf Setting and general assumptions on the coefficients}
We fix now an Hilbert space $H$, real and separable, we are going to study the following Lyapunov equation:
\begin{equation}\label{Lyap}
\left\{  \begin{array}{l}   -dP_t= (A' P_t + P_t A +C_t'Q_t + Q_t C_t  + C'_t P_t C_t) \,dt + S_t \, dt  - Q_t d\, W_t, \\ \\
\quad P_T=M
  \end{array}
\right .
\end{equation}
in the space $L(H)$.

%
%

The following assumptions on $A$, $C$, $S$ and $M$ will be used
throughout the paper:

\begin{hypothesis}\label{genhyp}

 $ $
\begin{enumerate}
\item[A1)]  $A$ is a self adjoint operator in $H$ and there exist a  complete orthonormal basis $\{ e_k: k \geq 1\}$ in $H$ (that we fix from now on), a sequence of real numbers $\{ \lambda_k: k\geq 1\}$ and $\omega \in \R$, such that
\begin{align} \label{diagonal}
A e_k  =- \lambda_k e_k, \quad \text{ with}
\quad\omega \leq \lambda_1  \leq \lambda_2 \leq \dots \leq   \lambda_k \leq \dots ,  \end{align}

\noindent Moreover we assume that for a suitable $\rho \in \big(\frac{1}{4},\frac{1}{2}\big)$, it holds
\begin{equation} \label{convautov}\sum_{k \geq 1}\lambda_k^{-2\rho} < +\infty.\end{equation}

\noindent Without weakening the generality of the problem we can, and will, assume that $\omega>0$ (just multiply $P$ and $Q$ by an exponential weight) .

\noindent As it is well known in this case $A$ generates an analytic
semigroup $(e^{tA})_{t \geq 0}$ with
$
|e^{tA}|_{L(H)} \leq 1 $.
\item[A2)]
We assume that $C\in
L^{\infty}_{\mathcal{P},S}(\Omega\times[0,T];L(H)).$ We denote
with $M_C$ a positive constant such that:
\begin{equation*}
|C(t,\omega)|_{L(H)} < M_C, \quad \mathbb{P}-\text{a.s. and for
a.e. } t \in (0,T).
\end{equation*}
\item[A3)]   $S \in L^\infty_{\P,S}( (0,T)\times \Omega;\Sigma^+(H)))$
 and $M \in L^\infty_S(\Omega,\F_T;\Sigma^+(H))$.
\end{enumerate}
\end{hypothesis}
$ $

\begin{remark}
We notice that requirement A.1) in \ref{genhyp}  is easily fulfilled in the case when
 $A$ is the realization of  the Laplace operator in $H= L^2([0,\pi])$ with Dirichlet boundary conditions.
One has indeed:
\begin{gather*}
D(A)= H^2 ([0,\pi]) \cap H^1_0([0,\pi]), \\
e_k (x)= (2/\pi)^{1/2}\sin{kx} , \quad  k=1, 2, \dots , \\
|\nabla e_k (x)| \leq  (2/\pi)^{1/2}  {k },   \quad    k=1, 2, \dots , \\
\lambda_k = k^2, \quad  k=1, 2, \dots .
\end{gather*}
Similar considerations can be done for the Laplace operator with Dirichlet boundary conditions on bounded domains of $\R^n$.

While requirement A.2) is fulfilled, for instance, as soon as  $C(t,\omega)$ is defined on $L^2([0,\pi])$ by $(C(t,\omega) x)(\xi):= c (t,\omega, \xi)x(\xi)$, with $c$ any bounded and progressive measurable map
  $[0,T ]\times \Omega \times [0,\pi]\rightarrow \mathbb{R}$.  The same holds for A.3), see also section 10 of \cite{GuaTess}.

\end{remark}
\noindent{\bf The  Hilbertian triple $V \hookrightarrow_d H  \hookrightarrow_d V'$ }

In this paragraph we introduce the Hilbertian triple we will use to build the effective Hilbert space of operators where we are going to solve the Lyapunov equation.
Let
\begin{equation}\label{defK}
V := D((-A)^\rho)=
\{ x \in H:  \sum_{n = 1}^\infty \lambda_n^{2\rho} |\langle x, e_n\rangle | ^2:=|x|^2_V<\infty\}.
\end{equation}
By construction $V$ is an Hilbert space endowed with its natural scalar product, in par\-ti\-cu\-lar $ \{ \lambda_n^{-\rho} e_n\}_{ n\geq 1}$ is a complete orthormal basis in $V$.

We can consider  also its topological dual $K'$ that has the following characterization:
\begin{equation}\label{defK'}
V' := D((-A)^{-\rho})
\end{equation}
Notice that $V'$ is the completion of $H$ with the norm
$ | \cdot |^2_{V'}= \sum_{n = 1}^\infty \lambda_n^{-2\rho} |\langle x, e_n\rangle | ^2$ and $ \{ \lambda_n^{\rho} e_n\}_{ n\geq 1}$ and that is a complete orthormal basis in $V'$.

Once we make the usual identification $ H  \simeq H' $,  we have the following dense inclusions:
\begin{equation}\label{tripletta}
 V \hookrightarrow_d H \hookrightarrow_d V'
\end{equation}
We notice that both inclusion operators are  Hilbert-Schmidt class

\begin{remark}
Under the previous hypotheses \ref{genhyp}, we have for all $t>0$
\begin{gather}\label{regsem}
 t^{\rho}|e^{tA}|_{L(H,V)} \leq 1,\;\;   t^{\rho}|e^{tA}|_{L(V',H)} \leq M_A, \\ \label{propsem} |e^{tA}|_{L(V)} \leq 1,\;\;  |e^{tA}|_{L(V')} \leq 1 .
\end{gather}

\end{remark}

\noindent{\bf The  Hilbert space $\K$. }
We set
\begin{equation}
\K := L_2(V;H) \cap L_2(H;V'),
\end{equation}
where   $L_2(V;H) $ denotes the Hilbert space of Hilbert-Schmidt operators
form $V$ to $H$, endowed with the Hilbert-Schmidt norm
$|T|_{L_2(V;H)}= (\sum_{i=1}^{\infty} |T f_i|^2_H)$
($\{f_i:i\in\mathbb{N}\}$ being a complete orthonormal basis-b.o.c.-in $V$), see \cite{DPZ1}.
$\K$ will be endowed with the natural norm $|T|_{\K}^2=|T|_{L_2(V;H)}^2+|T|_{L_2(V;H)}^2$

The obvious similar definition holds for $L_2 (H;V')$. \\
At last we introduce the following subspace of $K$:
\begin{equation}
\K_s := \{ G \in  L_2(V;H) \cap L_2(H;V')\text{ such that } \langle G x, y\rangle_H=\langle x, G y\rangle_H \text{ for all } x,y \in V\}
\end{equation}
We resume its main properties in the following Lemma.
\begin{lemma}\label{propK}
The following hold:
\begin{enumerate}
\item[(i)] $\mathcal{K}$ is a separable Hilbert space,
\item[(ii)]  $L(H) \subset  \mathcal{K}$,
\item [(iii)] $T \in \mathcal{K}$ iff $ T \in L(V;H) \cap  L(H;V')$ and $|T|_{\K}^2 = \sum_{k =1}^{\infty} \lambda_k ^{-2\rho} (|Te_k|^2 _H + |T'e_k|^2_H) < \infty$, where $T'\in L(V;H) \cap  L(H;V')$ is the adjoint of $T$ (in the sense that $\<Tv, w\>=\<v, T'w\>.$
    whenever $v\in V$ and
    $w\in H$ or $w\in V$ and
    $v\in H$.)
    \item [(iv)] If $T \in \K_s$ then $|T|_{\K_s}^2 = 2\sum_{k =1}^{\infty} \lambda_k ^{-2\rho} |Te_k|^2 _H $
\end{enumerate}
\end{lemma}
\Dim We omit the proof of $(i)$, being obvious.

$(ii)$ Let $G \in L(H)$, then  since $ \{ \lambda_n^{-\rho} e_n\}_{ n\geq 1}$ is a basis of $V$, we have:
\begin{equation}
|G |_{L_2(V;H)}=\Big( \sum_{n=1}^\infty \lambda_n^{-2\rho} |G e_n|_H^2 \Big)^{1/2} \leq |G|_{L(H)} \Big(\sum_{n=1}^\infty \lambda_n^{-2\rho} \Big)^{1/2}
\end{equation}
Moreover, recalling that $ \{ e_n: n \geq 1\}$ is a b.o.c. of $H$, we have:
\begin{align}
& \nonumber|G |_{L_2(H;V')}=\Big( \sum_{n=1}^\infty |G e_n|_{V'}^2 \Big)^{1/2} \leq |G|_{L(H)} \Big(\sum_{n=1}^\infty\sum_{h=1}^\infty \lambda_h^{-2\rho}|\langle e_n, e_h\rangle | ^2 _H\Big)^{1/2}\\ & =|G|_{L(H)} \Big(\sum_{h=1}^\infty\lambda_h^{-2\rho} \sum_{n=1}^\infty|\langle e_n, e_h\rangle | ^2 _H\Big)^{1/2} = |G|_{L(H)} \Big(\sum_{h=1}^\infty\lambda_h^{-2\rho}  \Big)^{1/2}
\end{align}
Thus $G \in \K$.

$(iii)$ Notice that, for any b.o.c. $\{ f_k: k \geq 1\}$ of $H$, we have:
\begin{align}
\sum_{k=1}^\infty |Tf_k|^2_{V'} = \sum_{k=1}^\infty  \sum_{h=1}^\infty \lambda_h^{-2\rho}\langle f_k,T' e_h \rangle _H^2=
 \sum_{h=1}^\infty\sum_{k=1}^\infty  \lambda_h^{-2\rho}\langle f_k, T'e_h \rangle _H^2= \sum_{h =1}^{\infty} \lambda_h ^{-2\rho}  |T'e_h|^2_H.
\end{align}
\finedim
\section{Mild Solutions of the Lyapunov Equation}

The natural space in which the deterministic Lyapunov equation is
studied is the space $\Sigma(H)$ of bounded self adjoint operators
in $H$. Unfortunately this is not an Hilbert space and this  fact causes
serious difficulties when considering stochastic backward
differential equations  (for instance the essential tool
given by the Martingale Representation Theorem does not hold). To
overcome this difficulty we will work in the bigger space $\K$ that is a separable Hilbert space.

For convenience we rewrite the equation of interest:
\begin{equation}\label{Lyap1}
\left\{  \begin{array}{l}   -dP_t= (A' P_t + P_t A +C'Q_t + Q_t C  + C' P_t C) \,dt + S_t \, dt  - Q_t d\, W_t, \\ \\
\quad P_T=M
  \end{array}
\right .
\end{equation}
\begin{definition}\label{defmild}
A {\em mild solution} of problem \eqref{Lyap1} is a couple of
processes
$$(P,Q) \in L^2_{\mathcal{P},S}(\Omega,C([0,T];\Sigma(H)))
\times L^2_\mathcal{P}(\Omega\times [0,T];\K_s)$$
that solves the following equation, for all $t \in [0,T]$:
\begin{align}\label{Lyapmild}
P(t)&=e^{(T-t)A'}Me^{(T-t)A}+ \int_t^Te^{(s-t)A'}
S(s)e^{(s-t)A}\,ds\nonumber \\
&+ \int_t^Te^{(s-t)A'}
\Big[C'(s)P(s)C(s)+C'(s)Q(s)+Q(s)C(s)\Big]e^{(s-t)A}\,ds
\\\nonumber
 &+ \int_t^Te^{(s-t)A^*} Q(s)e^{(s-t)A}\,dW(s)\quad\quad \mathbb{P}-\text{a.s.}
\end{align}
\end{definition}
We first prove an a-priori estimate for mild solutions.
\begin{proposition}\label{unicitaloc}
Let $(P,Q)$ a mild solution to \eqref{Lyapmild}.
Then there exists a $\delta_0 >0$ just depending on $T$ and the constants $M_A, M_C$ and $\rho$ introduced in \ref{genhyp} such that for every $0 \leq \delta \leq \delta_0$ the following holds:
\begin{align}\label{apriori}
|P|^2_{L^2(\Omega;C([T-\delta,T]; L(H)))}+ \E\int_{T-\delta}^T|Q(s)|^2_{\K}\,ds \leq c \Big(\E|M|^2_{L(H)} + \delta \E\int_{T-\delta}^T |S(s)|^2 _{L(H)}\,ds\Big).
\end{align}
where $c$ is a positive constant depending on $\delta_0, M_A, M_C, \rho$ and $T$.
\end{proposition}
\Dim
Let $(P,Q) \in L^2_{\mathcal{P},S}(\Omega,C([0,T];L(H)))
\times L^2_\mathcal{P}(\Omega\times [0,T];\K_s)$ be any mild solution, hence we have that:
\begin{align}\label{Lyapmildcond}
P(t)&=\Et \Big[e^{(T-t)A}Me^{(T-t)A} +\int_t^Te^{(s-t)A}
S(s)e^{(s-t)A}\,ds\Big] \\\nonumber
&+\Et\Big[ \int_t^Te^{(s-t)A}
\Big(C'(s)P(s)C(s)+C'(s)Q(s)+Q(s)C(s)\Big)e^{(s-t)A}\,ds\Big]
\quad\quad \mathbb{P}-\text{a.s.}
\end{align}
We notice that if $(L(t))_{T\geq 0}$ is a Banach space valued process then by Doob's $L^2$ inequality
$$ \E \sup_{t \in [r,T]}|\Et L(t)|^2\leq   \E \sup_{t \in [r,T]}[\Et( \sup_{t \in [r,T]} |L(t)|)]^2
\leq4 \E \sup_{t \in [r,T]}|L(t)|^2$$
Moreover we have:
\begin{align}
& \E \sup_{t \in [r,T]}|e^{(T-t)A}Me^{(T-t)A} |^2_{L(H)} \leq  |M|^2_{L(H)} \\
& \E \sup_{t \in [r,T]}\Big| \int_t^Te^{(s-t)A}C'(s)P(s)C(s) e^{(s-t)A}\,ds \Big|^2_{L(H)} \leq
M_C^4(T-r) \E\int_r^T |P(u)|_{L(H)}^2\,ds  \\
& \E \sup_{t \in [r,T]}\Big| \int_t^Te^{(s-t)A}S(s) e^{(s-t)A}\,ds \Big|^2_{L(H)}
\leq (T-r)E\int_{r}^T |S(s)|_{L(H)}^2 \,ds \label{stimaS}
\end{align}
In estimating the latter terms we notice that even if $G \in \K$ it is not true in general that $GC \in \K$,
 therefore we have to use the regularity properties of the semigroup \eqref{regsem}.
\begin{align*}
&  \E \sup_{t \in [r,T]} \Big|\int_t^Te^{(s-t)A} \Big[C'(s)Q(s)+Q(s)C(s)\Big]e^{(s-t)A}\,ds\Big|^2_{L(H)} \leq \\
&  2 \E \Big\{ \sup_{t \in [r,T]} \Big[ \int_t^T |e^{(s-t)A} C'(s)Q(s)e^{(s-t)A}|_{L(H)} \,ds \Big]^2+\sup_{t \in [r,T]} \Big[ \int_t^T |e^{(s-t)A} Q(s)C(s)e^{(s-t)A}|_{L(H)} \,ds \Big]^2
\Big\}
\end{align*}
Let us consider the first term:
\begin{align*}
& \E \Big\{ \sup_{t \in [r,T]} \Big[ \int_t^T |e^{(s-t)A} C'(s)Q(s)e^{(s-t)A}|_{L(H)} \,ds \Big]^2  \\&\leq
\E \sup_{t \in [r,T]}\Big [ \int_t^T|e^{(s-t)A}|_{L(H)} |C'(s)|_{L(H)} |Q(s)|_{L(V,H)} |e^{(s-t)A}|_{L(V)}\, ds \Big ]^2 \\ &\leq M_C^2(T-r)\E\int_r^T |Q(s)|^2_\K \, ds
\end{align*}
and the second one:
\begin{align}\label{stimaaprioriQ}
& \E \sup_{t \in [r,T]} \Big[ \int_t^T |e^{(s-t)A} Q(s)C(s)e^{(s-t)A}|_{L(H)} \,ds \Big]^2 \nonumber  \\ &\leq
\E \sup_{t \in [r,T]}\Big [ \int_t^T|e^{(s-t)A}|_{L(V';H)} |Q(s)|_{L(H;V')}|C(s)|_{L(H)} |e^{(s-t)A}|_{L(H)}\, ds \Big ]^2 \nonumber   \\&\leq  \E  \sup_{t \in [r,T]}\Big(\int_t^T\frac{M_C}{(s-t)^\rho} |Q(s)|_{\K}\, ds\Big)^2 \leq M_C^2(T-r) ^{1-2\rho} \int_r^T |Q(s)|_{\K}^2 \, ds.
\end{align}
Summing up all these estimates we obtain that, for $r=T-\delta$:
\begin{align} \label{StimaPUnif_zero}
& \E \sup_{u\in[T-\delta,T]} |P(u)|_{L(H)}^2
\\&\nonumber \leq C \Big( |M|^2_{L(H)} +\delta^2\E \sup_{u\in[T-\delta,T]} |P(u)|_{L(H)}^2 + \delta^{1-2\rho} \E \int_{T-\delta}^T |Q(s)|_{\K}^2 \,ds  +\delta \E \int_{T-\delta}^T |S(s)|_{L(H)}^2 \,ds\Big)
\end{align}
where $C$ depends only on $M_C, \rho$ and $T$ and for $\delta$ small enough (changing the value of the constant $C$)
\begin{equation} \label{StimaPUnif}
\E \sup_{u\in[T-\delta,T]} |P(u)|_{L(H)}^2
\leq C \Big( |M|^2_{L(H)} + \delta^{1-2\rho} \E \int_{T-\delta}^T |Q(s)|_{\K}^2 \,ds  +\delta \E \int_{T-\delta}^T |S(s)|_{L(H)}^2 \,ds\Big)
\end{equation}
Now we have to recover an estimate for $Q$, this can not be done in the same way because the term $Q(s)C(s) \notin \K$, and we can not follow the technique introduced in \cite{HuPeng1991}.

Therefore we exploit some duality relation. First of all we multiply both sides by the linear operators  $J_n:=n (nI-A)^{-1}$.

Such family of operators have the following properties:

\begin{enumerate}
\item $J_n e_k= \frac{n}{(n+\lambda_k)}e_k, \text{  for every } k\geq 1, \quad n \geq 1$,\\
\item $|J_n|_{L(H)} \leq 1, \quad |J_n|_{L(V)} \leq 1, \quad |J_n|_{L(V')} \leq 1,$ for every $n \geq 1$,\\
\item$ |J_n|_{L(H,V)} \leq   n^{\rho}, \quad |J_n|_{L(V',H)} \leq  n^{\rho}, $\\
\item $\lim_{n \to +\infty} J_n x=x, \text{ for every } x \in H$, \\
\item $J_n \in L_2(H)$, for every $n \geq 1$, and $ |J_n|_{L_2(H)} \leq |I_{{V,H}}|_{L_2(H)} |J_n|_{L(H,V)}$.\\
\end{enumerate}
hence equation \eqref{Lyapmild}, setting $P^n(s)= J_nP(s)J_n$ and $Q^n(s)= J_nQ(s)J_n$ becomes:
\begin{align}\label{mildappr}
\nonumber P^n(t)&=e^{(T-t)A} J_nMJ_ne^{(T-t)A}
+ \int_t^Te^{(s-t)A}
J_nC'(s)P(s)C(s)J_n e^{(s-t)A}\,ds\\\nonumber&  + \int_t^Te^{(s-t)A}
J_nS(s)J_ne^{(s-t)A}\,ds +  \int_t^Te^{(s-t)A}\Big [J_nC'(s)Q(s)J_n+J_nQ(s)C(s)J_n\Big]e^{(s-t)A}\,ds
\\
 &+ \int_t^Te^{(s-t)A} Q^n(s)e^{(s-t)A}\,dW(s)\quad\quad \mathbb{P}-\text{a.s.}
\end{align}
Notice that, thanks to the regularization property of $J_n$, $(P^n,Q^n) \in  L^2_\mathcal{P}(\Omega\times [0,T];L_2(H))\times L^2_\mathcal{P}(\Omega\times [0,T];L_2(H))$.
In particular
\begin{align*}
|Q_n(s)|^2_{L_2(H)}\leq  |J_n|^2_{L(V';H)} |Q(s)|^2_\K
\end{align*}
Moreover $(P^n,Q^n)$ is also the unique mild solution of:
\begin{equation}\label{Lyap1n}
\left\{  \begin{array}{l}   -dP^n_t= (A' P^n_t + P^n_t A +C'Q^n_t + Q^n_t C  + C' P^n_t C) \,dt + \hat{S}^n_t \, dt  - Q^n_t d\, W_t, \\ \\
\quad P_T=M^n
  \end{array}
\right .
\end{equation}
where $\hat{S}^n_s=J_nC'_sP_sC_sJ_n+J_nS_sJ_n +J_nC'_sQ_sJ_n+J_nQ_sC_sJ_n \in  L^2_\mathcal{P}(\Omega\times [0,T];L_2(H))$. We wish to apply Lemma $2.1$  of \cite{HuPeng1991}.
Let us check that $ \hat{S}^n$ has the required $L_2$ regularity:
\begin{multline}
\E \int_0^T |J_nC'(s)P(s)C(s)J_n|^2_{ L_2(H)} \,ds \leq \E \int_0^T |J_n|^2_{L(H)} |C'(s)|^2_{L(H)}  |P(s) |^2_{L(H)}|C(s)|^2_{L(H)} |J_n|^2_{L_2(H)}
\,ds  \\ \leq M_C^4 |J_n|^2_{L_2(H)}|P|_{L^2(\Omega;C([T-\delta,T];L(H)))}
\end{multline}
\begin{multline}
\E \int_0^T |J_nQ(s)C(s)J_n|^2_{ L_2(H)} \,ds \leq \E \int_0^T |J_n|^2_{L(V',H)}  |Q(s) |^2_{L_2(H;V')}|C(s)|^2_{L(H)} |J_n|^2_{L_2(H)}
\,ds  \\ \leq n^{2\rho}M_C^2 \E  \int_0^T|Q(s)|^2_{\K}\, ds.
\end{multline}
\begin{multline}
\E \int_0^T |J_nC'(s)Q(s)J_n|^2_{ L_2(H)} \,ds \leq \E \int_0^T    |J_n|^2_{L(H)}|C'(s)|^2_
{L(H)}|Q(s) |^2_{L_2(V;H)}  |J_n|^2_{L(H;V)}\,ds  \\ \leq n^{2\rho}M_C^2 \E  \int_0^T|Q(s)|^2_{\K}\, ds.
\end{multline}
We seek for an estimate independent of $n$ for the martingale term. We are going to use a duality argument, with this purpose we introduce an operator valued process defined as follows
\begin{equation}\label{defLn}
L^n(s) e_k := 2\lambda_k^{-2\rho} Q^n(s) e_k, \qquad \text{ for } k \geq 1.
\end{equation}
Let us fix $\delta >0$ then consider the following process
\begin{equation}\label{defXn}
X^n_t = \int_{T-\delta}^t e^{(t-s)A} L^n (s)e^{(t-s)A} \, dW(s), \qquad  t \in [T-\delta,T].
\end{equation}
It  can be  easily verified that $ X^n \in  C_\P([T-\delta,T];L^2(\Omega; L_2(H)))$.
Therefore, by standard regularization arguments, see for instance \cite{DPZ1}  for the forward equation and   \cite{GuaTess}  for backward equation we can prove that:
\begin{align}\label{dualita}
&\E\langle  X^n(T), P^n (T)\rangle_{L_2(H)}=\E \int_{T-\delta} ^T  \langle L^n(s), Q^n(s) \rangle_{L_2(H)} \, ds - \E \int_{T-\delta} ^T  \langle X^n(s), J_nS(s)J_n \rangle_{L_2(H)}\,ds
\nonumber\\ & - \E \int_{T-\delta} ^T  \langle X^n(s), J_nC'(s)P(s)C(s)J_n+ J_nC'(s)Q(s)J_n+J_nQ(s)C(s)J_n \rangle_{L_2(H)}\,ds.
\end{align}
First of all notice that $\langle L^n(s),  Q^n(s)\rangle _{L_2(H)} = 2 \sum_{k =1}^{\infty} \lambda_k^{-2\rho} | Q^n(s) e_k|^2_H$,  such quantity corresponds to $|Q^n|^2_{\K}$ being $Q^n$ a symmetric operator.
Thus
\begin{equation}\label{stima3}
|\E \int_{T-\delta} ^T  \langle L^n(s), Q^n(s) \rangle_{L_2(H)} \, ds| = \E \int_{T-\delta}^T|Q^n(s) |^2_\K \, ds
\end{equation}
Let us estimate  the process $ X^n_T$, we have for every $ t \in [T-\delta,T]$:
\begin{align}\label{stimaXn}
\E \sum_{k\geq 1} |  X^n(t) e_k|^2_H\lambda_k^{2\rho}  &= \sum_{k \geq 1} \E \Big|\int_{T-\delta}^t e^{(t-s)A} L^n (s) e^{(t-s)A} e_k\, dW_s\Big|^2_H \lambda_k^{2\rho} \\\nonumber &=   \sum_{k \geq 1} \E \int_{T-\delta}^t  \lambda_k^{2\rho} | e^{(t-s)A} L^n (s) e^{(t-s)A}e_k |^2_H\,ds \\&\nonumber \leq  \E \int _{T-\delta}^T \sum_{k \geq 1} \lambda_k ^{-2\rho} 2 |Q^n(s) e_k|^2 _H\, ds =  \E \int _{T-\delta}^T  |Q^n(s)|^2 _\K \,ds
\end{align}
Therefore, using $\eqref{stimaXn}$ with $r= T-\delta$ we have
\begin{align}\label{stima2}& |\E\langle  X^n(T), P^n (T)\rangle_{L_2(H)}| \\&= |\E \sum_{k=1}^{\infty} \langle  X^n(T) e_k, P^n (T)e_k \rangle|\leq \Big( \E \sum_{k=1}^{\infty} | X^n(T) e_k|^2 \lambda_k^{2\rho}\Big)^{\frac{1}{2}} \Big( \E \sum_{k=1}^{\infty} | P^n(T) e_k|^2\lambda_ k^{-2\rho}\Big)^{\frac{1}{2}} \nonumber\\ &\leq
 C\Big( \E \int _{T-\delta}^T  |Q^n_s|^2 _\K \,ds \Big)^{\frac{1}{2}}\Big(\E  | P^n(T)|^2_{L(H)} \Big)^{\frac{1}{2}}\nonumber\
\end{align}
Moreover,  thanks to $\eqref{StimaPUnif}$ and $|P^n(T)|_{L(H)}\leq |P(T)|_{L(H)}$,  we end up with
\begin{align}\label{stima2-bis}& |\E\langle  X^n(T), P^n (T)\rangle_{L_2(H)}| \\ &\leq
 C\Big( \E \int _{T-\delta}^T  |Q^n_s|^2 _\K \,ds \Big)^{\frac{1}{2}}\Big( |M|^2_{L(H)}
 + \delta \E \int _{T-\delta}^T |S(s)|^2 ds+ \delta^{1-2\rho} \E \int_{T-\delta}^T |Q(s)|_{\K}^2 \, ds \Big)^{\frac{1}{2}}\nonumber\
\end{align}
Regarding $\displaystyle \E \int_{T-\delta} ^T  \langle X^n(s), J_n C'(s)Q(s)J_n+J_nQ(s)C(s)J_n \rangle_{L_2(H)}\,ds $ we have:
\begin{align*}
&\Big|  \E \int_{T-\delta} ^T  \langle X^n(s), J_n C'(s)Q(s)J_n\rangle_{L_2(H)}\,ds \Big|\leq
  M_C^2\E \int_{T-\delta} ^T\Big (\sum_{k\geq 1} | X^n(s) e_k|^2_H \lambda_k^{2\rho}\Big)^{\frac{1}{2}} |Q(s)|_\K \,ds \\
& \leq C  \E \int_{T-\delta} ^T\Big(  \int_{T-\delta} ^T  |Q^n(s)|^2_\K \,ds \Big)^{\frac{1}{2}} |Q(s)|_\K \,ds \leq \frac{1}{4} \E \int_{T-\delta} ^T  |Q^n(s)|^2_\K \,ds + C \delta \E \int_{T-\delta} ^T  |Q(s)|^2_\K \,ds
\end{align*}
with $C>0$ a constant that may change form line to line but always depends only on the ones introduced in \ref{genhyp}.
Notice that
\begin{align}\label{extraccia}
  & \E \int_{T-\delta} ^T \langle X^n(s),  J_nQ(s)C(s)J_n\rangle_{L_2(H)}  \, ds=  \E \int_{T-\delta} ^T  \sum_{k=1}^{\infty}  \langle X^n(s) e_k,  J_nQ(s)C(s)J_n e_k\rangle _H
\,ds \nonumber\\& =  \E \int_{T-\delta} ^T  \sum_{k=1}^{\infty}\sum_{h=1}^{\infty}   \langle e_k, X^n(s) e_h \rangle \langle e_k,  J_nC'(s)Q(s)J_n e_h\rangle _H \,ds\\& \leq \E \int_{T-\delta} ^T\sum_{h=1}^{\infty}
 |X^n(s) e_h| |J_nC'(s)Q(s)J_ne_h|\,ds  \nonumber\\ & \leq   \E \int_{T-\delta} ^T (\sum_{h=1}^{\infty}  \lambda_h^{2\rho}  |X^n(s) e_h|^2)^{1/2} (\sum_{h=1}^{\infty}  \lambda_h^{-2\rho}  |Q(s) e_h|^2)^{1/2} \, ds. \nonumber
\end{align}
Thus the same conclusion holds, so we have that, by \eqref{stimaXn}:
\begin{align}\label{stima1}
&\Big|  \E \int_{T-\delta} ^T \langle X^n(s),J_n C'(s)Q(s)J_n+ J_nQ(s)C(s)J_n\rangle_{L_2(H)}  \, ds \Big| \nonumber\\ & \leq \frac{1}{2} \E \int_{T-\delta} ^T  |Q^n(s)|^2_\K \,ds + C \delta \E \int_{T-\delta} ^T  |Q(s)|^2_\K \,ds
\end{align}
Moreover we have that
\begin{align}\label{stima4}
&\Big|  \E \int_{T-\delta} ^T \langle X^n(s), J_nC'(s)P(s)C(s)J_n\rangle_{L_2(H)}  \, ds \Big| \nonumber\\ & \leq C \delta  |P|_{L^2_\mathcal{P} (\Omega; C([T-\delta,T];L(H))) }^2+ \frac{1}{8}\E \int_{T-\delta}^T  |Q^n(s)|^2_\K  \,ds,
\end{align}
and that, similarly,
\begin{align}\label{stima5}
&\Big|  \E \int_{T-\delta} ^T \langle X^n(s), J_nS(s)J_n\rangle_{L_2(H)}  \, ds \Big|\leq C \E \int_{T-\delta}^T  |S(s)|^2_{L(H)} ds+ \frac{1}{8}\E \int_{T-\delta}^T  |Q^n(s)|^2_\K  \,ds,
\end{align}

Taking into account \eqref{stima2-bis}, \eqref{stima1}, \eqref{stima4}  and \eqref{stima5} we have that there exists a positive constant
$C$ independent of $n$ and $\delta$ such that
\begin{multline}
\label{stimaQnunif}
\E \int_{T-\delta} ^T  |Q^n(s)|^2_\K \,ds \leq C\Big( |M|^2_{L(H)}+ \delta \E \int_{T-\delta}^T  |S(s)|^2_{L(H)} ds +\delta^{1-2\rho}\E \int_{T-\delta}^T |Q(s)|_{\K}^2 \, ds \Big)
\end{multline}
From \eqref{StimaPUnif} and \eqref{stimaQnunif} the claim follows  since $|Q^n(s)|^2_\K \nearrow |Q(s)|^2_\K$ choosing a $\delta$ small enough such that $C\delta^{1-2\rho} < \frac{1}{2}$. \finedim

With identical argument we get the estimate in the easier case in which the term $C'PC$ is not present
\begin{remark}\label{stima-a-priori-P-noto}
Assume that $Q \in L^2_\mathcal{P}(\Omega\times [0,T];\K_s)$ and that $P$ given by
\begin{align}\label{Lyapmild-no-CPC}
P(t)&=e^{(T-t)A'}Me^{(T-t)A}+ \int_t^Te^{(s-t)A'}
S(s)e^{(s-t)A}\,ds\nonumber \\
&+ \int_t^Te^{(s-t)A'}
\Big[C'(s)Q(s)+Q(s)C(s)\Big]e^{(s-t)A}\,ds
+ \int_t^Te^{(s-t)A^*} Q(s)e^{(s-t)A}\,dW(s)\quad\quad \mathbb{P}-\text{a.s.}
\end{align}
is an adapted $\K$-valued process.

 Then there exists a $\delta_0 >0$ just depending on $T$ and the constants $M_C$ and $\rho$ introduced in \ref{genhyp} such that for every $0 \leq \delta \leq \delta_0$ the following holds:
\begin{align}\label{apriori-no-CPC}
|P|^2_{L^2(\Omega;C([T-\delta,T]; L(H)))}+ \E\int_{T-\delta}^T|Q(s)|^2_{\K}\,ds \leq c \Big(\E|M|^2_{L(H)} + \delta \E\int_{T-\delta}^T |S(s)|^2 _{L(H)}\,ds\Big).
\end{align}
with $c$ is a positive constant depending on $\delta_0, M_C, \rho$ and $T$.
\end{remark}
We are now in a position to prove existence and uniqueness of the solution to the mild Lyapunov equation
\begin{theorem}\label{lyapunov.teo}
Under  assumptions \ref{genhyp}  equation
\eqref{Lyap1} has a unique mild solution $(P,Q)$.
\end{theorem}
\Dim
The idea is classical: we will buid a map  $\Gamma$ from the space  $L^2_\mathcal{P}(\Omega,C([0,T];H))$ into its self and prove that is a contraction for small time.

In completing this program we follow three steps.

{\bf Step 1: regularization}
We introduce some regularizing processes in order to define $\hat{P}=\Gamma(P)$ for an arbitrary $P\in L^2_\mathcal{P}(\Omega,C([0,T];\Sigma(H)))$.
So we fix $P$ and for every $n \geq 1$ we consider the following problem: find $\hat{P}^n,\hat{Q}^n$ such that
\begin{align}\label{Lyapmildappr}
\nonumber \hat{P}^n(t)&=e^{(T-t)A} J_nMJ_ne^{(T-t)A}
+ \int_t^Te^{(s-t)A}
C'(s) J_nP(s)J_nC(s) e^{(s-t)A}\,ds\\\nonumber&  + \int_t^Te^{(s-t)A}
J_nS(s)J_ne^{(s-t)A}\,ds +  \int_t^Te^{(s-t)A}(C'(s) \hat{Q}^n(s)+ \hat{Q}^n(s)C(s))e^{(s-t)A}\,ds
\\
 &+ \int_t^Te^{(s-t)A} \hat{Q}^n(s)e^{(s-t)A}\,dW(s)\quad\quad \mathbb{P}-\text{a.s.}.
\end{align}
Notice that for every $n \in \mathbb{N}$, we have that $C' J_nPJ_nC,  J_nSJ_n  \in L^2_\mathcal{P}(\Omega\times [0,T];L_2(H))$,  $J_nMJ_n \in L_2(H)$.
Moreover for every $C \in L(H)$, the map $G \in L_2(H) \to C' G +GC \in L_2(H)$ is Lipschitz continuous.

Thus Lemma $2.1$  of \cite{HuPeng1991} applies and we can deduce that there exists
a unique solution $(\hat{P}^n,\hat{Q}^n) \in  L^2_\mathcal{P}(\Omega\times [0,T];L_2(H))\times
L^2_\mathcal{P}(\Omega\times [0,T];L_2(H))$ to eq.\eqref{Lyapmildappr}.
Moreover by Remark \ref{stima-a-priori-P-noto} there exists $\delta_0 <1 $ small enough and  independent of $n$ such that $\forall \delta \leq \delta_0$
\begin{align} \label{StimaPQnUnif}
& \E\!\!\!\!  \sup_{u\in[T-\delta,T]}\!\!\!\! |\hat{P}^n(u)|_{L(H)}^2+ \E \int_{T-\delta} ^T\!\!\!\! |\hat{Q}^n(s)|^2_\K \,ds \
\ \leq C \Big( |M|^2_{L(H)} +\delta^2\E \sup_{u\in[r,T]} |P(u)|_{L(H)}^2 + \delta \int_{r}^T |S(s)|_{L(H)}^2 \,ds \Big),
\end{align}
with $C$ a constant depending only on $M_C, T$ and $\rho$ but not on $n$.

\smallskip

We notice here that the operator $P \rightarrow C'PC$ is lipschitz from $L_2(H)$ to itself as well. We can not treat it as
the term  $G   \to C' G +GC $ since we will then need to lower the regularity of $P$ to the space $\K$
and if $P$ only belonges to $\K$ then the operator $e^{sA} C'PC e^{sA}$ is not well defined while
$G\rightarrow e^{sA}[ C'G+ GC ]e^{sA}$ is well defined from $\K$ to itself.

 {\bf Step 2: limiting procedure}
Let us evaluate the difference $\hat{P}^n-\hat{P}^m$ for two integers $m,n$:
\begin{align}\label{Lyapmildapprdiff}
&\nonumber \hat{P}^n(t) - \hat{P}^m(t)=e^{(T-t)A} J_nMJ_ne^{(T-t)A}- e^{(T-t)A} J_mMJ_me^{(T-t)A} \\&+ \int_t^Te^{(s-t)A}
(J_nS(s)J_n - J_mS(s)J_m)e^{(s-t)A}\,ds\nonumber
\\& \nonumber + \int_t^Te^{(s-t)A}
C'(s)(J_nP(s)J_n- J_mP(s)J_m)C(s)e^{(s-t)A}\, ds\\\nonumber &+
 \int_t^Te^{(s-t)A}\Big[C'(s)(\hat{Q}^n(s)- \hat{Q}^m(s))+(\hat{Q}^n(s)- \hat{Q}^m(s))C(s)\Big]e^{(s-t)A}\,ds
\\
 &+ \int_t^Te^{(s-t)A} [\hat{Q}^n(s)-\hat{Q}^m(s)]e^{(s-t)A}\,dW(s)\quad\quad \mathbb{P}-\text{a.s.}
\end{align}
we are going to show that
\begin{align}\label{limPn}
&\lim_{m,n \to \infty }\E\!\sup_{t\in [T-\delta, T]}  |\hat{P}^n(t)-\hat{P}^m(t)|^2 _\K= 0 \\\label{limQn}
&\lim_{m,n \to \infty } \E \int_{T-\delta}^{T} |\hat{Q}^n(s)-\hat{Q}^m(s)|^2 _\K\, ds=0
\end{align}
Let' s begin to prove \eqref{limPn}
by noticing that:
\begin{align}\label{Lyapmildattesacond}
&\nonumber \hat{P}^n(t) -\hat{P}^m(t)=\Et (e^{(T-t)A} J_nMJ_ne^{(T-t)A}- e^{(T-t)A} J_mMJ_me^{(T-t)A} )\\&+ \Et \Big(\int_t^Te^{(s-t)A}
(J_nS(s)J_n - J_mS(s)J_m)e^{(s-t)A}\,ds \Big)\nonumber
\\& \nonumber +\Et \Big( \int_t^Te^{(s-t)A}
C'(s)(J_nP(s)J_n- J_mP(s)J_m)C(s)e^{(s-t)A}\, ds \Big) \\\nonumber &+
\Et\Big( \int_t^Te^{(s-t)A}\Big[C'(s)(\hat{Q}^n(s)- \hat{Q}^m(s))+(\hat{Q}^n(s)- \hat{Q}^m(s))C(s)\Big]e^{(s-t)A}\,ds\Big),
\quad \mathbb{P}-\text{a.s.}
\end{align}
Being $M$ a symmetric operator, we have that
\begin{align*}
|e^{(T-t)A} (J_nMJ_n-  J_mMJ_m)e^{(T-t)A} |^2_\K =  \sum_{k =1}^{\infty} \lambda_k ^{-2\rho}|e^{(T-t)A} (J_nMJ_n-  J_mMJ_m)e^{(T-t)A} e_k|^2_H
\end{align*}
For every fixed $k \geq 1$:
\begin{equation*}\label{M1}
\lim_{n,m \to \infty}| (J_nM(J_n-  J_m) e_k|^2_H=0,\quad \forall t \in [0,T], \qquad \mathbb{P} -\text{a.s.}
\end{equation*}
and \begin{equation*}\label{M2}
\lim_{n,m \to \infty}|(J_n-  J_m)M J_m e_k|^2_H=0, \quad \forall t \in [0,T] ,\qquad \mathbb{P} -\text{a.s.}
\end{equation*}
Moreover
 \begin{equation*}\label{M3}
 \sum_{k =1}^{\infty} \lambda_k ^{-2\rho}|(J_nMJ_n-  J_mMJ_m)e_k|^2_H \leq  M_A^4M_M^2\sum_{k =1}^{\infty} \lambda_k ^{-2\rho}  < \infty, \qquad \mathbb{P} -\text{a.s.}.
\end{equation*}
Hence by the Dominated Convergence Theorem and the Doob inequality for martingales:
 \begin{align}\label{MFin}
&\lim_{n,m \to \infty} {\E}[\sup_{t \in [T-\delta,T]}|\Et(e^{(T-t)A}(J_nMJ_n-  J_mMJ_m)e^{(T-t)A})|^2_\K] \\\nonumber& \leq
\lim_{n,m \to \infty}4 \,\E |(J_nMJ_n-  J_mMJ_m) |^2_\K =0.
\end{align}
The second and the third term are similar so we'll give the details only of the third.
As before we have that for every $k\geq 1$:
\begin{equation*}
\lim_{n,m \to \infty}|(C'(s)(J_nP(s)J_n- J_mP(s)J_m)C(s)e_k|^2_H=0, \quad  \mathbb{P} -\text{a.s.}  \text{ and for  a.e. } s\in [T-\delta,T]
\end{equation*}
and   $\mathbb{P}$ -{a.s.} and for  a.e. $ s\in [T-\delta,T]$,
\begin{align*}
&
\sum_{k\geq 1} \lambda_k^{-2\rho}|(C'(s)(J_nP(s)J_n - J_mP(s)J_m)C(s)e_k|^2_H\,ds\leq  M_C^4\sum_{k\geq 1} \lambda_k^{-2\rho} < \infty. \end{align*}
Therefore again by the Dominated Convergence Theorem  and the Doob inequality for martingales:
 \begin{align}\label{CPCFin}
&\lim_{n,m \to \infty} {\rm E}\sup_{t\in[T-\delta,T]}\Big| \Et\Big(\int_t^Te^{(s-t)A}
[(C'(s)(J_nP(s)J_n- J_mP(s)J_m)C(s)]e^{(s-t)A}\,ds\Big)\Big|^2_\K\nonumber  \\&  \leq \delta \lim_{n,m \to \infty} 4 \, \E \int_{T-\delta}^T
\sum_{k\geq 1} \lambda_k^{-2\rho}|(C'(s)(J_nP(s)J_n - J_mP(s)J_m)C(s)e_k|^2_H\,ds =0.
\end{align}
At last  let us consider the term
\begin{align*}& \E \sup_{t \in [T-\delta,T]}\Big|\Et \Big( \int_t^Te^{(s-t)A}[ C'(s)(\hat{Q}^n(s)-\hat{Q}^m(s))+(\hat{Q}^n(s)-\hat{Q}^m(s))C(s)]e^{(s-t)A}\, ds\Big) \Big|_\K^2
\end{align*}
First of all:
\begin{align*}
&\Big( \int_t^T|e^{(s-t)A}(\hat{Q}^n(s)-\hat{Q}^m(s))C(s)e^{(s-t)A}|_{\K}\, ds \Big)^2
\\ & =\Big[ \int_t^T(\sum_{k \geq 1}\lambda_k^{-2\rho}|e^{(s-t)A}(\hat{Q}^n(s)-\hat{Q}^m(s))C(s)e^{(s-t)A}e_k|^2_H)^{1/2}\, ds \Big]^2 \\ & \leq
\Big(\int_t^T|e^{(s-t)A}|_{L(V';H)}|(\hat{Q}^n(s)-\hat{Q}^m(s))|_{L_2(H;V')} (\sum_{k \geq 1}\lambda_k^{-2\rho}|C(s)e^{(s-t)A}e_k|^2_H)^{1/2}\, ds\Big)^2  \\ &\leq M_C^2 ( \sum_{k \geq 1} \lambda_k^{-2\rho}) \int_t
^T(s-t)^{-2\rho}\, ds \int_t^T |\hat{Q}^n(s)-\hat{Q}^m(s)|^2_{\K} \, ds \\&\leq
C \delta^{1-2\rho} \int_{T-\delta}^T   |\hat{Q}^n(s)-\hat{Q}^m(s)|^2_{\K} \, ds\end{align*}
Similarily
\begin{align*} &\Big[ \int_t^T(\sum_{k \geq 1}\lambda_k^{-2\rho}|e^{(s-t)A}C'(s)(\hat{Q}^n(s)-\hat{Q}^m(s))e^{(s-t)A}e_k|^2_H)^{1/2}\, ds \Big]^2 \\ & \leq
\Big(\int_t^T|e^{(s-t)A}|_{L(H)} |C'(s)|_{L(H)}^2 (\sum_{k \geq 1}\lambda_k^{-2\rho}e^{-2(s-t)\lambda_k}|(\hat{Q}^n(s)-\hat{Q}^m(s))e_k|^2_H)^{1/2}\, ds\Big)^2  \\ &\leq M_C^2\delta \int_t^T |\hat{Q}^n(s)-\hat{Q}^m(s)|^2_{\K} \, ds \\&\leq
C \delta^{1-2\rho} \int_{T-\delta}^T   |\hat{Q}^n(s)-\hat{Q}^m(s)|^2_{\K} \, ds\end{align*}
Hence:
\begin{align}\label{stimahatPmn}
&\E\!\sup_{t\in [T-\delta, T]}  |\hat{P}^n(t)-\hat{P}^m(t)|^2 _\K \leq C\Big[ \delta^{1-2\rho} \int_{T-\delta}^T \E  |\hat{Q}^n(s)-\hat{Q}^m(s)|^2_{\K} \, ds + \E |(J_nMJ_n-  J_mMJ_m) |^2_\K\nonumber  \\&  +
 \delta  \, \E \int_{T-\delta}^T
|C'(s)(J_nP(s)J_n - J_mP(s)J_m)C(s)|^2_\K\,ds +\E \int_{T-\delta}^T
|J_nS(s)J_n - J_mS(s)J_m)|^2_\K\,ds  \Big] \nonumber\\& \leq C \delta^{1-2\rho} \E\int_{T-\delta}^T |\hat{Q}^n(s)-\hat{Q}^m(s)|^2_{\K} \, ds + R(m,n)
\end{align}
with $R(m,n) \to 0$ as $m,n \to +\infty$.

The duality relation between $\hat{P}^n-\hat{P}^m$ and $\hat{X}^n-\hat{X}^m$ yields to:
\begin{align}\label{dualitamn}
&\nonumber\E\langle  \hat{X}^n(T) - \hat{X}^m(T), \hat{P}^n (T) - \hat{P}^m(T)\rangle_{L_2(H)}=\E \int_{T-\delta} ^T  \langle \hat{L}^n(s) -\hat{L}^m(s), \hat{Q}^n(s) -\hat{Q}^m(s)\rangle_{L_2(H)} \, ds \\\nonumber &- \E \int_{T-\delta} ^T  \langle \hat{X}^n(s) -\hat{ X}^m(s), J_nS(s)J_n -J_mS(s)J_m\rangle_{L_2(H)}\,ds
\nonumber\\ & - \E \int_{T-\delta} ^T  \langle \hat{X}^n(s) -\hat{X}^m(s),  C'(s)J_nP(s)J_nC(s)-C'(s)J_mP(s)J_mC(s)\rangle_{L_2(H)}\,ds\nonumber \\& -
\E \int_{T-\delta} ^T  \langle X^n(s) -X^m(s),
 C'(s)(\hat{Q}^n(s)-\hat{Q}^m(s))+(\hat{Q}^n(s)-\hat{Q}^m(s))C(s)\rangle_{L_2(H)}\,ds.
\end{align}
where
$$\E \int_{T-\delta} ^T  \langle \hat{L}^n(s) -\hat{L}^m(s), \hat{Q}^n(s) -\hat{Q}^m(s)\rangle_{L_2(H)} \, ds = \E \int_{T-\delta} ^T| \hat{Q}^n(s) -\hat{Q}^m(s) |^2_\K\,ds.$$
As in  \eqref{stimaXn} we have:
\begin{align}\label{stimaXnm}
\E \sum_{k\geq 1} |  (\hat{X}^n(t)-  \hat{X}^m(t))e_k|^2_H\lambda_k^{2\rho}  \leq  \E \int _{T-\delta}^T  |\hat{Q}^n(s)-\hat{Q}^m(s)|^2 _\K \,ds,
\end{align}
where $\hat{X}^n$ and $\hat{X}^m$ are defined as in \eqref{defXn} with $Q_n$ replaced by $\hat{Q}^n$ and
we get, noticing that $|\langle X,Z \rangle_{L_2(H)}| \leq (\sum_{k=1}^\infty |X e_k|^2 \lambda_k^{2\rho})^{1/2} |Z|_\K$
\begin{align}
&\E \int_{T-\delta} ^T  |\hat{Q}^n(s)-\hat{Q}^m(s)|^2_\K \,ds \leq \E\!\sup_{t\in [T-\delta, T]} \! |\langle \hat{X}^n(T)- \hat{X}^m(T),\hat{P}^n(T)-\hat{P}^m(T) \rangle_{L_2(H)}| \nonumber \\ &
+\E \int_{T-\delta}^T | \langle J_nS(s)J_n -J_mS(s)J_m, \hat{X}^n(s)- \hat{X}^m(s)\rangle _{L_2(H)}|\,ds
\nonumber
\\&+ \E\int_{T-\delta}^T |\langle C'(s)[J_nP(s)J_n-J_mP(s)J_m]C(s), \hat{X}^n(s)- \hat{X}^m(s)\rangle_{L_2(H)} |\,ds  \\\nonumber &+\E\int_{T-\delta}^T |\langle X^n(s) -X^m(s),
C'(s)(\hat{Q}^n(s)-\hat{Q}^m(s))\rangle_{L_2(H)} | \,ds
\\  \nonumber &+  \E\int_{T-\delta}^T|\langle X^n(s) -X^m(s),(\hat{Q}^n(s)-\hat{Q}^m(s))C(s)\rangle _{L_2(H)}|\,ds
 \Big)
= I_1+ I_2+ I_3+ I_4+I_5.
\end{align}
We have
\begin{gather*}
I_1 \leq \E (\sum_{k=1}^\infty \lambda_k^{2\rho} |(\hat{X}^n(T)- \hat{X}^m(T))e_k|^2)^{1/2} |\hat{P}^n(T)-\hat{P}^m(T)|_\K  \\ \leq \frac{l}{2}
  \E \int _{T-\delta}^T  |\hat{Q}^n(s)-\hat{Q}^m(s)|^2 _\K \,ds + \frac{1}{2l}|\hat{P}^n(T)-\hat{P}^m(T)|^2_\K
\end{gather*}
\begin{gather*}
I_2+I_3  \leq  \frac{l}{2}\E \int _{T-\delta}^T \sum_{k=1}^\infty \lambda_k^{2\rho} |(\hat{X}^n(s)- \hat{X}^m(s))e_k|^2 \, ds + \frac{1}{2l} \E \int _{T-\delta}^T| \langle J_nS(s)J_n -J_mS(s)J_m|_\K^2\, ds\\
 +  \frac{1}{2l}\E \int _{T-\delta}^T|C'(s)[J_nP(s)J_n-J_mP(s)J_m]C(s)|^2_\K \, ds,
\end{gather*}
\begin{gather*}
I_4 \leq \frac{1}{2l}\E \int _{T-\delta}^T \sum_{k=1}^\infty \lambda_k^{2\rho} |(\hat{X}^n(s)- \hat{X}^m(s))e_k|^2 \, ds +\frac{l}{2} \E \int _{T-\delta}^T  |\hat{Q}^n(s)-\hat{Q}^m(s)|^2 _\K \,ds \\
\leq  \frac{\delta}{2l} \E \int _{T-\delta}^T  |\hat{Q}^n(s)-\hat{Q}^m(s)|^2 _\K \,ds +\frac{l}{2} \E \int _{T-\delta}^T  |\hat{Q}^n(s)-\hat{Q}^m(s)|^2 _\K \,ds
\end{gather*}
$I_5$ can be treated as $I_4$, following \eqref{extraccia}.

Summarizing and choosing $l$ small enough (depending only on the constants introduced in \ref{genhyp}), we finally get
\begin{align}\label{stimaQnmunif}
&\E \int _{T-\delta}^T  |\hat{Q}^n(s)-\hat{Q}^m(s)|^2 _\K \,ds  \leq C\Big( \E |\hat{P}^n(T)-\hat{P}^m(T)|^2 _\K
+ \E \int_{T-\delta}^T |J_nS(s)J_n -J_mS(s)J_m|^2_\K\,ds
\\&\nonumber+ \E\int_{T-\delta}^T |C'(s)J_nP(s)J_nC(s)-C'(s)J_mP(s)J_mC(s)|^2_\K\,ds
+ \delta \, \E\int_{T-\delta}^T |\hat{Q}^n(s)-\hat{Q}^m(s)|^2 _\K \, ds \Big)
\end{align}
Putting together \eqref{stimahatPmn} and \eqref{stimaQnmunif} we then prove that for a small enough $\delta$:
\begin{gather*}
\lim_{m,n \to \infty }\E\!\sup_{t\in [T-\delta, T]}  |\hat{P}^n(t)-\hat{P}^m(t)|^2 _\K= 0  \\
\lim_{m,n \to \infty } \E \int_{T-\delta}^{T} |\hat{Q}^n(s)-\hat{Q}^m(s)|^2 _\K \, ds=0
.\end{gather*}
Therefore there exist the limit $\hat{P} \in L^2_\mathcal{P} (\Omega; C([T-\delta,T];\K)) $ and $\hat{Q} \in L^2_\mathcal{P} (\Omega\times [T-\delta,T];\K))$ such that:
\begin{gather}
\lim_{n \to \infty }\E\!\sup_{t\in [T-\delta, T]}  |\hat{P}^n(t)-\hat{P}(t)|^2 _\K= 0  \label{hatPlim}\\
\lim_{m, \to \infty } \E \int_{T-\delta}^{T} |\hat{Q}^n(s)-\hat{Q}(s)|^2 _\K \, ds=0\label{hatQlim}
.\end{gather}

{\bf Step 3: construction of $\Gamma$}.
Being the equation linear, thanks to \eqref{hatPlim} and \eqref{hatQlim}, we obtain the following relation:
\begin{align}\label{Gamma}
\nonumber \hat{P}(t)&=e^{(T-t)A}Me^{(T-t)A}
+ \int_t^Te^{(s-t)A}
C'(s)P(s)C(s) e^{(s-t)A}\,ds\\\nonumber&  + \int_t^Te^{(s-t)A}
S(s)e^{(s-t)A}\,ds +  \int_t^Te^{(s-t)A}(C'(s)\hat{Q}(s)+\hat{Q}(s)C(s))e^{(s-t)A}\,ds
\\
 &+ \int_t^Te^{(s-t)A} \hat{Q}(s)e^{(s-t)A}\,dW(s)\quad\quad \mathbb{P}-\text{a.s.}
\end{align}
The fact that $\hat{P} \in  L^2_{\mathcal{P},S} (\Omega; C([T-\delta,T];L(H)))$ follows from Remark \ref{stima-a-priori-P-noto}.
So far we have that the map $\Gamma$ such that $\Gamma(P)= \hat{P}$ is actually defined from the space $ L^2_{\mathcal{P},S} (\Omega; C([T-\delta,T];L(H)) $ into itself.

{\bf Step 4: $\Gamma$ is a contraction for a suitable $\delta$}.
Let $P^1$ and $ P^2$ two elements of  $ L^2_{\mathcal{P},S} (\Omega; C([T-\delta,T];L(H))$, then we can evaluate the difference between $\Gamma(P^1)$ and $\Gamma(P^2)$. Indeed we have:
\begin{align}\label{Differenza}
\nonumber (\hat{P}^1- \hat{P}^2)(t)&=\int_t^Te^{(s-t)A}
C'(s)(P^1-P^2)(s)C(s) e^{(s-t)A}\,ds\\\nonumber&  + \int_t^Te^{(s-t)A}[C'(s)(\hat{Q}^1-\hat{Q}^2)(s)+(\hat{Q}^1-\hat{Q}^2)(s)C(s)]e^{(s-t)A}\,ds
\\
 &+ \int_t^Te^{(s-t)A} (\hat{Q}^1-\hat{Q}^2)(s)e^{(s-t)A}\,dW(s)\quad\quad \mathbb{P}-\text{a.s.}
\end{align}
Clearly \eqref{StimaPUnif} and \eqref{stimaQnunif} hold also in this case
\begin{align}\label{StimaP12Unif}
& \E \sup_{u\in[T-\delta,T]} |(\bar{P}^1- \bar{P}^2)(u)|_{L(H)}^2
\\&\nonumber \leq C \Big(\delta\, \E \sup_{u\in[T-\delta,T]} |(P^1-P^2)(u)|_{L(H)}^2 + \delta^{1-2\rho} (\E \int_{T-\delta}^T |(\hat{Q}^1-\hat{Q}^2)(u)|_{\K}^2 \,du \Big),
\end{align}
with the constant $C$ depending on the constants $ M_C$ and $T$ but not on $\delta$.
And the same holds for $\hat{Q}^1-\hat{Q}^2$:
\begin{multline}
\label{stimaQ12unif}
\E \int_{T-\delta} ^T  |(\hat{Q}^1-\hat{Q}^2)(s)|^2_\K \,ds \leq C\Big( \delta
 |P^1-P^2|_{L^2_\mathcal{P} (\Omega; C([T-\delta,T];L(H))) }^2  +\delta^{1-2\rho} \E \int_{T-\delta}^T | (\hat{Q}^1-\hat{Q}^2)(s)|_{\K}^2 \, ds \Big)
\end{multline}
So we can find a $\delta$ small enough such that $\Gamma$ is a contraction and there's a fixed point $P$ . The couple $(P,\hat{Q})$, where $\hat{Q}$ is defined in  \eqref{Gamma} is the mild solution in $[ T-\delta,  T]$.

\medskip
{\bf Step 5: construction of the mild solution}
Since the problem  is linear  and the value of $\delta$ depends only on the constants introduced in \ref{genhyp}, can restart on $[T-2\delta, T-\delta]$ with final datum $P(T-\delta)$. Proceeding backwards we arrive to cover the whole interval $[0,T]$.

\medskip
{\bf Step 6: uniqueness}
From Proposition \ref{unicitaloc} we have that there is local uniqueness for the mild solution.
Being $\delta_0$ independent of the data,  we can deduce global uniqueness.
\finedim

We end the section proving the following stability results for the approximants processes  $\hat{P}^n$:
\begin{proposition}\label{stabilita}
Under the hypotheses of the previous theorem,  let $\hat{P}^n$ defined by \eqref{Lyapmild} and $P$ the mild solution just obtained, then the following holds
there exists a $\delta >0$ such that  for every $\varepsilon < \delta_1$:
\begin{align}\label{limiteLH}
\lim_{n \to \infty}\E \sup_{t \in [T-\delta,T-\varepsilon]} |P(t)- \hat{P}^n(t)|^2_{L(H)}=0.
\end{align}
\end{proposition}
\Dim
For every $t \in [0,T]$ we have
\begin{align}\label{LyapmildPPn}
&P(t)- \hat{P}^n(t)= \Et \big\{ e^{(T-t)A'}(M- J_nMJ_n)e^{(T-t)A}+ \int_t^Te^{(s-t)A'}
(S(s)-J_nS(s)J_n)e^{(s-t)A}\,ds +\\ \nonumber
&\int_t^Te^{(s-t)A'}
[C'(s)(P(s)-J_n{P}(s)J_n)C(s)+C'(s)(Q(s)-\hat{Q}^n(s))+(Q(s)-\hat{Q}^n(s))C(s)]e^{(s-t)A}\,ds\big\},
\end{align}
thus, assume that $\delta <1$
\begin{align*}
& \E \sup_{t \in [T-\delta,T-\varepsilon]}| \Et e^{(T-t)A}(M- J_nMJ_n)e^{(T-t)A} |^2_{L(H)}\\
 &  = \E \sup_{t \in [T-\delta,T-\varepsilon]}| \Et e^{(T-\varepsilon - t)A}e^{\varepsilon A}(M- J_nMJ_n)e^{\varepsilon A} e^{(T-\varepsilon - t)A}|^2_{L(H)}
 \\& \leq 4  \varepsilon ^{- 2\rho} \, \E |M- J_nMJ_n|^2_{\K}, \\
& \E \sup_{t \in [T-\delta,T - \varepsilon]}\Big|\Et \int_t^{T}e^{(s-t)A}
C'(s)(P(s)- J_nP(s)J_n)C(s) e^{(s-t)A}\,ds \Big|^2_{L(H)} \\&\leq
4\delta ^{1-2\rho} \E \int_{T-\delta}^T|C'(s)( P(s)-J_n{P(s)}J_n) C(s)|_\K^2\,ds, \\
&  \E \sup_{t \in [T-\delta,T-\varepsilon]} \Big|\Et \int_t^Te^{(s-t)A} [ C'(s)(Q(s) -\hat{Q}^n(s))+(Q(s)-\hat{Q}^n(s))C(s)]e^{(s-t)A}\,ds \Big|^2_{L(H)}\\ &\leq 2\  \E  \sup_{t \in [T-\delta,T]}\Big(\Et \int_t^T\frac{M_C}{(s-t)^\rho} |(Q(s) -\bar{Q}^n(s)|_{\K}\, ds\Big)^2 \leq
8 M_C^2 \, \delta ^{1-2\rho} \E \int_{T-\delta}^T |(Q(s) -\bar{Q}^n(s)|_{\K}^2 \, ds, \\
& \E \sup_{t \in [T-\delta,T]}\Big| \int_t^Te^{(s-t)A}(J_nS(s)J_n- S(s))e^{(s-t)A}\,ds \Big|^2_{L(H)}
\leq \delta^{1-2\rho}\E\int_{r}^T |S(s)- J_nS(s)J_n|_{\K}^2 \,ds.
\end{align*}
Summing up all these estimates we deduce that there exists a constant $C$ depending only on $M_C, \rho $ such that :
\begin{align}\label{stimaLH}
& \E \sup_{t \in [T-\delta,T-\varepsilon]} |P(t)- \hat{P}^n(t)|^2_{L(H) }\leq  C( \varepsilon ^{- 2\rho}  \E |M- J_nMJ_n|^2_{\K} + \delta ^{1-2\rho} \E \int_{T-\delta}^T| P(s)-J_n{P(s)}J_n|_\K^2\,ds\nonumber \\&
 \delta ^{1-2\rho} \E \int_{T-\delta}^T |Q(s) -\hat{Q}^n(s)|_{\K}^2 \, ds+ \delta^{1-2\rho}\E\int_{r}^T |S(s)- J_nS(s)J_n|_{\K}^2 \,ds).
\end{align}

Thanks to previous considerations in particular \eqref{hatQlim}, and recalling that by dominated convergence theorem
$ \E \int_{T-\delta}^T| P(s)-J_n{P(s)}J_n|_\K^2 \rightarrow 0$, we deduce the thesis.
\finedim

\section{Backward Stochastic Riccati Equations and  LQ  Optimal Control}\label{Sec-LQ}
Besides hypotheses \ref{genhyp}, let us fix $T> S>0$ and consider the following infinite dimensional stochastic control problem, with {\em state equation} given  by
\begin{equation}\label{stato}
\left\{
\begin{array}{ll}
dy(t)=(Ay(t)+B(t)u(t)) \,dt + C(t)y(t)\,dW(t) &   S \leq r \leq t \leq T \\
y(r)=x,
\end{array}
\right.
\end{equation}
 where $u$ is the \emph{control} and takes values in another Hilbert space
 $U$.

Besides hypothesis \ref{genhyp} we assume that

\begin{enumerate}\label{hypB}
\item[A4)]
We assume that $B\in
L^{\infty}_{\mathcal{P},S}(\Omega\times[0,T];L(U;H)).$ We denote
with $M_B$ a positive constant such that:
\begin{equation*}
|B(t,\omega)|_{L(U;H)} < M_B, \quad \mathbb{P}-\text{a.s. and for
a.e. } t \in (0,T).
\end{equation*}
\end{enumerate}

We recall the definition of {\em mild solution}.
\begin{definition}\label{def-mild-stato}
 Given $x\in H$ and $u \in L^2_\P(\Omega
\times [t,T];U)$,
a mild solution of \eqref{stato} is a process $y
\in L^2_\P(\Omega \times [t,T];H)$ such that,
  almost everywhere in $\Omega \times [t,T]$,
\begin{align*}
y(s)= e^{(s-t)A}x+
\int_t^{s} e^{(s-\sigma)A}B(\sigma)u(\sigma)\,d\sigma +
\int_t^{s} e^{(s-\sigma)A}C(\sigma)y(\sigma)\,dW(\sigma).
\end{align*}
\end{definition}
The following existence and uniqueness 
 result holds:
\begin{theorem}\label{mildstato}
Assume \ref{genhyp}. Given any $x \in H$ and $u \in L^2_\P(\Omega
\times [t,T];U)$ 
 problem \eqref{stato} has a unique mild solution $y \in C_{\P}([t,T];L^2(\Omega;H))$.
Moreover,
\begin{equation}\label{stima-stato}
 \sup_{s \in [t,T]}\media |y(s)|^2 \leq C_2 \Big[|x|^2 +
 \media  \int_t^T|u(s)|^2\,ds \Big]
\end{equation}
for a suitable constant $C_2$ depending only on $T,M_B,M_C$ (notice that $C_2\geq 1$).

Finally if  $p >2$ and
$$ \media \Big( \int_t^T|u(s)|^2\,ds \Big)^{\frac{p}{2}} < \infty,$$
then we have that $y\in L^p_{\P}(\Omega;C([t,T];H))$ and
\begin{equation}\label{stima-stato-p}
 \media \sup_{s \in [t,T]}|y(s)|^p \leq C_p\bigg[|x|^p +
 \media \bigg( \int_t^T|u(s)|^2\,ds \bigg)^{\frac{p}{2}}\bigg]
\end{equation}
for some positive constant $C_p$ depending on $p,T,M_B,M_C$.
\end{theorem}
The {\em cost functional}  to minimize over all processes taking values in $L^2_\P(\Omega \times [0,T], U)$- the  space of {\em admissible controls} is
\begin{equation}\label{costo}
\media \int_0^T \left(|\sqrt{S}(s)y(s)|_H^2 + |u(s)|_H^2\right)\,ds +
\media \langle My(T),y(T)\rangle_H.
\end{equation}
Associated to this {\em Linear and Quadratic} control problem we have the  following Backward Stochastic Riccati Equation (BSRE), see \cite{Bi76, Peng} and \cite{GuaTess} for the present infinite dimensional version:
\begin{equation}
 \label{Riccati}
\left\{
\begin{array}{l}
-dP(t)=(AP(t)+P(t)A+ C'(t)P(t)C(t)+C'(t)Q(t)+Q(t)C(t))\,dt    \\ \\
\qquad\qquad \;\;- (P(t)B(t)B^*(t)P(t)-S(t))\,dt  + Q(t)\,dW(t)  \qquad\qquad \qquad  t \in [0,T] \\ \\
P(T)=M
\end{array}
\right.
\end{equation}
In this section we will prove that such equation has a unique mild solution, in the sense of definition \ref{defmild}, improving the result obtained in \cite{GuaTess}. To be more specific we have
\begin{definition}\label{defmildRic}
A {\em mild solution} of problem \eqref{Riccati} is a couple of
processes
$$(P,Q) \in L^2_{\mathcal{P},S}(\Omega,C([0,T];\Sigma(H)))
\times L^2_\mathcal{P}(\Omega\times [0,T];\K_s)$$
that solves the following equation, for all $t \in [0,T]$:
\begin{align}\label{Ricmild}
P(t)&=e^{(T-t)A'}Me^{(T-t)A}+ \int_t^Te^{(s-t)A'}
S(s)e^{(s-t)A}\,ds\nonumber \\
&+ \int_t^Te^{(s-t)A'}
\Big[C'(s)P(s)C(s)-P)(s)B(s)B'(S)P(s)+C'(s)Q(s)+Q(s)C(s)\Big]e^{(s-t)A}\,ds
\\\nonumber
 &+ \int_t^Te^{(s-t)A^*} Q(s)e^{(s-t)A}\,dW(s)\quad\quad \mathbb{P}-\text{a.s.}
\end{align}
\end{definition}
We have indeed:
\begin{theorem}\label{main}
Assume that hypotheses \ref{genhyp}  and $\rm{A_4)}$ hold true.
Then there exists a unique mild solution $(P,Q)$ of equation \eqref{Riccati} in $ [0,T]$. Moreover $P \in  L^\infty_{\P,S}(\Omega\times (0,T);\Sigma^+(H))$.
Moreover,
 fix $T >0$ and $x\in H$, then
\begin{enumerate}
\item[{\rm 1.}] there exists a unique control $\overline u \in L^2_\P(\Omega\times[0,T];U)$ such that
\begin{equation*}
J(0,x, \overline u)= \inf _{u \in L^2_\P(\Omega\times[0,T];U) }
J(0,x,u);
\end{equation*}
\item[{\rm 2.}] if $\overline y$ is the mild solution of the state equation
corresponding to $\overline u$ (that is, the optimal state), then
$\overline y$ is the unique mild solution to the \emph{ closed
loop } equation 
\begin{equation}\label{loop.ban}
\left\{
\begin{array}{@{}ll}
d\overline y(r)=[A\overline y(r)-B(r)B'(r)P(r)\overline y(r)]\,dr+
 C\overline y(r)\,dW(r), \\
\overline y(0)=x;
\end{array}
\right.
\end{equation}
\item[{\rm 3.}] the following feedback law holds $\mathbb{P}$-a.s. for
almost every $s$:
\begin{equation}\label{feedback.ban}
     \overline u(s)=-B'(s)P(s)\overline y(s);
\end{equation}
\item[{\rm 4.}] the optimal cost is given by $ J(0,x,\overline u)= \langle P(0)x,x \rangle_H$.
\end{enumerate}
\end{theorem}
Before going into the details of the proof, we establish the following-priori estimate.
\begin{proposition}\label{stimapos}
Let $(\bar P,\bar Q )$ a mild solution of equation \eqref{Riccati} in $[\tau,T] \subset [0,T]$ such that $\bar P \in L^{\infty}_{\P,S}(\Omega \times [\tau,T], \Sigma(H)) $, then the following holds for every $ t \in [\tau,T]$:
\begin{itemize}
\item[(i)] for all $t \in [\tau,T]$,  $\bar{P}(t) \in \Sigma^+(H), \quad \mathbb{P} -a.s.$.
\item[(ii)]for all $t \in [\tau,T]$, \begin{equation}\label{stimaRiccati}
|\bar P(t)|_{L(H)} \leq  C_2( |M|_{L^\infty_{\P,S}(\Omega, \F_T; L(H))}+ (T-\tau) |S|_{ L^{\infty}_{\P,S}(\Omega \times [\tau,T], L(H))} )\qquad \mathbb{P}- \text{a.s.}
\end{equation}
where $C_2$ is given in \eqref{stima-stato}.
\end{itemize}

\end{proposition}
\Dim
{ \bf Step 1 [Fundamental relation for the Lyapunov equation]}.
Let $(P,Q)$ be the unique mild solution to the Lyapunov equation
\eqref{Lyapmild} and let $y^{t,x}$ be the mild solution to
\eqref{stato}, we claim that for all $t\in [0,T]$, $x \in H$,
it holds,
\begin{align}\label{rel-fond_0}
 \nonumber  \langle P(t)x,x \rangle _H &=
    \E^{\mathcal{F}_t}\langle  M y^{t,x}(T),y^{t,x}(T)\rangle+\E^{\mathcal{F}_t}\int_t^T\langle
    S(s)y^{t,x}(s),y^{t,x}(s)\rangle_H\,
    ds \\ &- 2  \, \Et \int_t^T \langle P(s) B'(s)y^{t,x}(s), u(s)\rangle \, ds , \hspace{3cm} \mathbb{P}\hbox{-}\text{\rm a.s.}.
\end{align}
Let us prove the claim.
We will use again the approximants processes  $(\hat{P}^n,\hat{Q}^n)$ introduced in the proof of theorem \ref{lyapunov.teo}. From
 proposition \ref{stabilita} we know that there's a $\delta$ small enough such that for every $\varepsilon < \delta$:
\begin{align}\label{limiteLH-bis}
\lim_{n \to \infty}\E \sup_{t \in [T-\delta,T-\varepsilon]} |P(t)- \hat{P}^n(t)|^2_{L(H)}=0.
\end{align}
On the other hand we have already noticed that $(\hat{P}^n,\hat{Q}^n)$ is a solution in the sense of Proposition 2.1 of
\cite{HuPeng1991}, therefore by Theorem 5.6 of \cite{GuaTess} we have that:
for all $t\in [0,T]$, $x\in H$, it holds,
$\mathbb{P}$-a.s., that
\begin{align}\label{rel-fondn}\nonumber
  \displaystyle \langle \hat{P}^n(t)x,x \rangle _H &=
    \E^{\mathcal{F}_t} \langle \hat{P}^n(T-\varepsilon)y^{t,x}(T-\varepsilon),y^{t,x}(T-\varepsilon)
    \rangle_H+\E^{\mathcal{F}_t}\int_t^{T-\varepsilon}\langle
    S(s)y^{t,x,u}(s),y^{t,x,u}(s) \rangle_H \, ds\\ \nonumber
&  +  \displaystyle \E^{\mathcal{F}_t}\int_t^{T-\varepsilon}\langle
 [ C'(s)\hat{P}^n(s)C(s) - C'(s)J_n{P}(s)J_n C(s)] y^{t,x,u}(s) ,y^{t,x,u}(s) \rangle _H\, ds  \\
&   \displaystyle - 2  \, \Et \int_t^{T-\varepsilon} \langle \hat{P}^n(s) B'(s)y^{t,x,u}(s), u(s)\rangle_H \, ds
 \end{align}
By (\ref{limiteLH-bis}) and  recalling that  $y\in L^p_{\P}(\Omega;C([t,T];H))$, $p\geq 2$ (see  (\ref{stima-stato-p})) we get that
$$\int_t^{T-\varepsilon} \langle
    \hat{P}^n(s)C(s)y^{t,x,u}(s),C(s)y^{t,x,u}(s) \rangle ds \rightarrow
    \int_t^{T-\varepsilon}\langle
    {P}(s)C(s)y^{t,x,u}(s),C(s)y^{t,x,u}(s) \rangle ds $$
    in $L^1$ norm. Moreover, since $\E \sup_{t \in [0,T]} |P(t)|^2_{L(H)} < +\infty$, by Dominated convergence theorem we obtain that
    $$\int_t^{T-\varepsilon} \langle
    {P}(s)J_nC(s)y^{t,x,u}(s),J_nC(s)y^{t,x,u}(s) \rangle ds \rightarrow
    \int_t^{T-\varepsilon} \langle
    {P}(s)C(s)y^{t,x,u}(s),C(s)y^{t,x,u}(s) \rangle ds $$
    again in $L^1$ norm.

Thus letting $n$ tend to $\infty$ in (\ref{rel-fondn}), we obtain that for
every $t \in [T-\delta, T]$, $\mathbb{P}$-a.s.:
\begin{align}\label{rel-fond}
\nonumber   \langle P(t)x,x \rangle _H &=
    \E^{\mathcal{F}_t} \langle {P}(T-\varepsilon)y^{t,x}(T-\varepsilon),y^{t,x}(T-\varepsilon) \rangle_H+\E^{\mathcal{F}_t}\int_t^{T-\varepsilon}\langle
    S(s)y^{t,x,u}(s),y^{t,x,u}(s) \rangle \, ds\\ & -2  \, \Et \int_t^{T-\varepsilon} \langle{P}(s) B'(s)y^{t,x,u}(s), u(s)\rangle \, ds
\end{align}
Now, thanks again  to $\E \sup_{t \in [0,T]} |P(t)|^2_{L(H)} < +\infty$, we can let $\varepsilon$  going to $0$
and get that for every $x \in H$, and every $ t \in [T-\delta,T]$, $\mathbb{P}$-a.s.:
\begin{align}\label{rel-fondfin}
 \nonumber  \langle P(t)x,x \rangle _H=
    \E^{\mathcal{F}_t} \langle My^{t,x}(T),y^{t,x}(T) \rangle_H+\E^{\mathcal{F}_t}\int_t^{T}\langle
    S(s)y^{t,x,u}(s),y^{t,x,u}(s) \rangle \, ds\\ -2  \, \Et \int_t^{T} \langle{P}(s) B'(s)y^{t,x,u}(s), u(s)\rangle \, ds
\end{align}
Choose $u=0$ then, see also Theorem 5.6 of \cite{GuaTess} we get that:
\begin{equation}\label{stimafondfin}
 \sup_{x \in H, \,  |x|_H=1} | \langle P(t)x,x \rangle _H | \leq C_2 ( |M|_{L^\infty_S(\Omega, \F_T,P)}
   + T |S| _{L^\infty_{\P,S}(\Omega\times (0,T);L(H))}), \quad \forall t \in  [T-\delta, T].\\
\end{equation}
We can prove relation \eqref{rel-fond} on the interval $[T-2\delta, T-\delta]$
(notice that $\delta$ does not depend on $M$) and so on to cover the  whole interval $[0,T]$, because $P(T-k\delta) \in L(H)$  for every $k=0,1,2,3, \dots$ and thus we can extend \eqref{stimafondfin} to the whole $[0,T]$.

\medskip
{\bf Step 2: upper bound}
Let $(\bar P, \bar Q)$ be the mild solution of the BSRE \eqref{Riccati} in $[\tau,T]$, we can see such couple of processes as the mild solution to the following Lyapunov equation, for  $t \in [\tau,T] $:
\begin{equation}
 \label{Lyapbar}
\left\{
\begin{array}{l}
-d\bar P(t)=(A\bar P(t)+\bar P(t)A+ C'(t)\bar P(t)C(t)+C'(t)\bar Q(t)+\bar Q(t)C(t) +\bar{S}(t))\,dt  +\bar Q(t)\,dW(t),\\ \\
\bar P(T)=M.
\end{array}
\right.
\end{equation}
with $\bar S = -B'\bar P \bar PB + S$, thus from \eqref{rel-fondfin}  and completing the square, we obtain
\begin{align}\label{relfondRiccati}
 \langle\bar P(t)x,x \rangle _H &=
    \E^{\mathcal{F}_t} \langle My^{t,x}(T),y^{t,x}(T) \rangle_H+\E^{\mathcal{F}_t}\int_t^{T}
    |u(s)|^2 \, ds \\  &+ \E^{\mathcal{F}_t}\int_t^{T}\langle
    S(s)y^{t,x,u}(s),y^{t,x,u}(s) \rangle \, ds - \Et \int_t^{T} |\bar{P}(s) B'(s)y^{t,x,u}(s)+u(s)|^2\, ds \nonumber
\end{align}
So, choosing the admissible control $u =0$, we get:
\begin{align}\label{relfondRiccati0}
 \langle \bar P(t)x,x \rangle _H &=
    \E^{\mathcal{F}_t} \langle My^{t,x,0}(T),y^{t,x,0}(T) \rangle_H+ \E^{\mathcal{F}_t}\int_t^{T}\langle
    S(s)y^{t,x,0}(s),y^{t,x,0}(s) \rangle \, ds \\  &- \Et \int_t^{T} |\bar{P}(s) B'(s)y^{t,x,0}(s)|^2\, ds \nonumber
\end{align}
From which we deduce the following upper bound
\begin{align}\label{Riccatiupper}
 \langle \bar P(t)x,x \rangle _H & \leq
   C_2 ( |M|_{L^\infty_S(\Omega, \F_T,P)}
   + T |S| _{L^\infty_{\P,S}(\Omega\times (0,T);L(H))}), \quad \forall t \in  [\tau, T].
\end{align}
{\bf Step 3: lower bound}
Let us consider the following equation for initial time $ t \in [\tau, T]$ and initial state $x$:
\begin{equation}\label{loopsmall}
\left\{
\begin{array}{ll}
d\overline y(s)=[A\overline y(s)-B(s)B'(s)\bar P(s)\overline y(s)]\,ds+
 C\overline y(s)\,dW(s),  \quad s \in [t,T] \\
\overline y(t)=x;
\end{array}
\right.
\end{equation}
Notice that,  thanks to the regularity of $\bar P$,  Theorem 3.2 of \cite{GuaTess} apply and in particular the following estimates holds true for the solution $\bar y^{t,x}$, for every $t \in [\tau, T]$:
\begin{equation}\label{stimastatoloop}
 \Et \sup_{s \in [t,T]}|\bar y(s)|^p \leq C_p |x|^p, \qquad \qquad \forall p \geq 2.
\end{equation}
where $C_p$ depends also on the $L^{\infty}$ norm of $\bar P$.
Therefore $ \bar u(s) = B'(s)\bar P(s)\bar{y}^{t,x}(s)$ is an admissible control, i.e. $\bar u \in L^2_{\P}(\Omega \times [t,T],U)$, and  \eqref{relfondRiccati} corresponds to
\begin{align}\label{relfondRiccatipos}
 \langle\bar P(t)x,x \rangle _H=  \E^{\mathcal{F}_t} \Big[\langle M\bar y^{t,x}((T),\bar y^{t,x}((T) \rangle_H\!+\!\int_t^{T}
  (  | B'(s)\bar P(s)\overline y^{t,x}((s)|^2 \!\!+\!|
   \sqrt{ S(s)}\bar y^{t,x}(s)|^2)\, ds\Big], \, \mathbb{P}-a.s.
\end{align}
Consequently from \eqref{relfondRiccatipos}  holding for every $t \in [\tau, T]$ we get $(i)$.
Eventually  \eqref{Riccatiupper} and \eqref{relfondRiccatipos} imply $(ii)$.
\finedim

We are now in the position to prove Theorem \ref{main}:

\medskip
{\bf Proof of Theorem \ref{main}.}

\noindent
{\bf Step 1: local existence and uniqueness} \\
In order to be able to follow the same argument not only on $[T-\delta,T]$
but also on  $[T-2\delta,T-\delta]$ and so on (with the same $\delta$) we prove existence of a solution (for notational convenience on $[T-\delta,T]$) with generic final condition $\tilde{M}\in L^\infty_{\P,S}(\Omega, \F_T; L(H))$
 with
 $$|\widetilde{M}|_{L^\infty_{\P,S}(\Omega, \F_T; L(H))}< C_2( |M|_{L^\infty_{\P,S}(\Omega, \F_T; L(H))}+ T|S|_{ L^{\infty}_{\P,S}(\Omega \times [0,T], L(H))} )$$
We fix  a number $ r $ with
$$r> C_2 ^2 |M|_{L^\infty_S(\Omega, \F_T,P)}
   + 2 C^2_2 T |S| _{L^\infty_{\P,S}(\Omega\times (0,T);L(H))}$$
   where $C_2$ is the the constant obtained in Proposition \ref{stimapos}
\begin{equation*}
\B(r)= \Big\{ P \in L^2 _{\P,S}(\Omega; C([T-\delta,T];L(H))):
 \sup_{t \in [T-\delta,T]}|P(t,\omega)|_{L(H)} \leq r\quad
\mathbb{P}\hbox{-}\text{\rm a.s.} \Big\}
\end{equation*}
where $ \delta >0$ will be fixed later on.
On $\B(r)$ we construct the  map $\Lambda: \B(r) \to \B(r)$, letting
$\Lambda(K)=P$, where $(P,Q)$ is the unique mild solution to
\eqref{Lyapmild} (in $[T- \delta,T]$) with $S$ replaced by $S-KBB^*K$
and $M$ by $\widetilde M$ that is verifies
\begin{eqnarray}
P(t)&=&e^{(T-t)A}\widetilde {M}e^{(T-t)A} + \int_t^Te^{(s-t)A}
[C'(s)P(s)C(s)+C'(s)Q(s)+Q(s)C(s)]e^{(s-t)A}\,ds\nonumber\\[1pt]
&&+\int_t^Te^{(s-t)A} S(s)e^{(s-t)A}\,ds+ \int_t^Te^{(s-t)A}K(s)B(s)B'(s)K(s) e^{(s-t)A}\,ds \nonumber\\[1pt]
&& +\int_t^Te^{(s-t)A} Q(s)e^{(s-t)A}\,dW(s) \nonumber\\[1pt]
 &&.\nonumber
\end{eqnarray} 
 First of all we check that it maps $\B(r)$ into itself. It is enough to show that for all $t\in [T-\delta,T]$ it holds
$|\Lambda(K)(t)|_{L(H)}\leq r$ $\mathbb{P}$-a.s. Thanks to
\eqref{stimaRiccati} we have that $\mathbb{P}$-a.s.
\begin{align*}
&|\Lambda( K )(t)|_{L(H)}\leq C_2 \bigg[|\widetilde M|_{L^{\infty}_S(\Omega,\F_T;L(H))}+
\delta |KBB'K|_{L^{\infty}_{\P,S}(\Omega \times [T-\delta,T]; L(H))}\nonumber\\
&+\delta |S|_{L^{\infty}_{\P,S}(\Omega \times [T-\delta,T]; L(H))})\,ds\bigg]\leq\\
& \leq  C_2^2 |M|_{L^{\infty}} + C_2 r^2 \delta M_B^2+ 2C_2^2T  |S|_{L^{\infty}_{\P,S}(\Omega \times [0,T]; L(H))} < r
\end{align*}
as soon as we choose $$\delta < \frac{r-(C_2^2 |M|_{L^{\infty}}+ 2C_2^2T  |S|_{L^{\infty}_{\P,S}(\Omega \times [0,T]; L(H))})}{C_2^2 M_B^2 r^2}.$$
Let $K_1$ and $K_2$ in $B(r)$, then by \eqref{rel-fond_0} evaluated at $ u=0$ we have:
\begin{align}\label{relfondiff}
 \langle  (P^1(t)-  P^2(t))x,x \rangle _H &=
   \E^{\mathcal{F}_t}\int_t^{T}\langle
    K^1(s)B(s)B'(s)(K^1(s) -K^2(s))y^{t,x,0}(s),y^{t,x,0}(s) \rangle \, ds \\  &-  \E^{\mathcal{F}_t}\int_t^{T}\langle
    K^2(s)B(s)B'(s)(K^1(s)-K^2(s))y^{t,x,0}(s),y^{t,x,0}(s) \rangle \, ds, \nonumber
\end{align}
thus, by H\"older inequality,
\begin{align}
&| \langle  (P^1(t)-  P^2(t))x,x \rangle _H| \leq
2   \E^{\mathcal{F}_t}\int_t^{T} r M_B^2 |K^1(s)-  K^2(s)|_{L(H)} |y^{t,x,0}(s)|^2 \,ds \\\nonumber &\leq
2 r M^2_B \int_t^{T} (\E^{\mathcal{F}_t}  |K^1(s)-  K^2(s)|^2_{L(H)})^{1/2}(\Et|y^{t,x,0}(s)|^4 )^{1/2}\,ds \\
\nonumber &\leq 2 r M^2_B \delta^2 (\sup_{t \in [T-\delta,T]}\E^{\mathcal{F}_t}  |K^1(s)-  K^2(s)|^2_{L(H)})^{1/2}
(\sup_{t \in [T-\delta,T]}\Et|y^{t,x,0}(s)|^4)^{1/2}
\end{align}
using again Doob inequality and \eqref{stimastatoloop} which we deduce:
\begin{equation}
\E \sup_{t \in [T-\delta, T]} |P^1(t)-  P^2(t) |^2 _{L(H)} \leq 16 r^2 M^4_B\delta^4 C_4 \E \sup_{t \in [T-\delta, T]} |K^1(t)-  K^2(t) |^2 _{L(H)}
\end{equation}
where $C_4= C_4 (r)$ is given in \eqref{stimastatoloop}.
Therefore reducing if necessary the value of $\delta$, we obtain that $\Lambda$ is a contraction.

{\bf Step 2: global existence and uniqueness.}
We notice that the choice of $\delta$ depends only on $r$ and the constants introduced in hypotheses \ref{genhyp}. Therefore we can repeat the previous step to cover the whole interval $[0,T]$.

{\bf Final step: synthesis of the optimal control.}
So far we have proved the existence and uniqueness of the mild solution for the BSRE, and thanks to Proposition \ref{stimapos} we also have that the first component of the solution $P \in L^{\infty}_{\P,S}(\Omega \times [0,T]; L(H))$. Consequently the closed loop equation \eqref{loop.ban} is well posed and the associated feedback control is admissible, hence the rest of the claims of the Theorem easily follow.

\finedim

\section{The Lyapunov Equation of the Maximum Principle}

In this section we extend Proposition \ref{unicitaloc} and Theorem \ref{lyapunov.teo} in order to cover the Lyapunov equation arising in the Maximum Principle for SPDE, see \cite{Fu-Hu-Te-1}, \cite[ Eq. (4.22)]{Fu-Hu-Te-2}.
We rewrite such equation with our notation
\begin{equation}\label{BSDEoperatorvalued}
    \left\{\begin{array}{l}
    -dP(t) = -Q(t) \,dW(t)
    +[
    A^*P(t) + P(t)A
     + A^*_{\sharp}(t)P(t) +P(t)A^*_{\sharp}]\,dt \\
    \qquad\qquad \, \,\, +[C(t)P(t) C(t)+ C(t)Q(t) +  Q(t)C(t)
    + S(t) ]\,dt
    \\
    P(T) =M,
\end{array}\right.
\end{equation}

where $A_{\sharp}\in
L^{\infty}_{\mathcal{P},S}((0,T)\times \Omega;L(H)).$ 

The presence of the bounded term $A_{\sharp}$ is completely irrelevant and we will not consider it in the following.

  On the contrary it is not 
possible, in this context, to require Assumption  $\rm{A3)}$.  Indeed Assumption  $\rm{A3)}$ has to be replaced by the weaker one 
\begin{hypothesis}
\begin{itemize}
\item[A3')]   $S \in L^2_{\P,S}( (0,T)\times \Omega;\K))$
 and $M \in L^\infty_S(\Omega,\F_T;L(H))$.
\end{itemize}
(notice that the assumption on $M$ remains unchanged)
\end{hypothesis}

Under ${\rm{A3')}}$ the statement of the a-priori estimate in Proposition \ref{unicitaloc} becomes:
\begin{proposition}\label{unicitalocLyap}
Let $(P,Q)$ a mild solution to \eqref{BSDEoperatorvalued}.
Then there exists a $\delta_0 >0$ just depending on $T$ and the constants $M_A,  M_C$ and $\rho$ introduced in $\rm{A1)-A2)}$ such that for every $0 \leq \delta \leq \delta_0$ the following holds:
\begin{align}\label{aprioriLyap}
|P|^2_{L^2(\Omega;C([T-\delta,T]; L(H)))}+ \E\int_{T-\delta}^T|Q(s)|^2_{\K}\,ds \leq c \Big(\E|M|^2_{L(H)} + \delta^{1-2\rho} \E\int_{T-\delta}^T |S(s)|^2 _{\K}\,ds\Big).
\end{align}
where $c$ is a positive constant depending on $\delta_0, M_A,  M_C, \rho$ and $T$.
\end{proposition}
\Dim
Let us reestimate \eqref{stimaS}. We have (by Cauchy inequality):
\begin{align}
 \E \sup_{t \in [r,T]}\Big| \int_t^Te^{(s-t)A}S(s) e^{(s-t)A}\,ds \Big|^2_{L(H)}
&\leq \E \sup_{t \in [r,T]} \Big(\int_t^T (s-t)^{-2\rho} \, ds  \int_{t}^T |S(s)|_{\K}^2 \,ds \Big) \\
&\leq (T-r)^{1-2\rho}\int_{r}^T |S(s)|_{\K}^2 \,ds \qquad\qquad \forall r \in [T-\delta,T] \nonumber\end{align}
Therefore \eqref{StimaPUnif_zero} becomes
\begin{align} \label{StimaPUnif_zeroLyap}
& \E \sup_{u\in[T-\delta,T]} |P(u)|_{L(H)}^2
\\&\nonumber \leq C \Big( |M|^2_{L(H)} +\delta^2\E \sup_{u\in[T-\delta,T]} |P(u)|_{L(H)}^2 + \delta^{1-2\rho} \E \int_{T-\delta}^T |Q(s)|_{\K}^2 \,ds  +\delta ^{1-2\rho}\E \int_{T-\delta}^T |S(s)|_{\K}^2 \,ds\Big)
\end{align}
From which we deduce:
\begin{align} \label{StimaPUnifLyap}
 \E \sup_{u\in[T-\delta,T]} |P(u)|_{L(H)}^2\leq C \Big( |M|^2_{L(H)}  + \delta^{1-2\rho} \E \int_{T-\delta}^T |Q(s)|_{\K}^2 \,ds  +\delta ^{1-2\rho}\E \int_{T-\delta}^T |S(s)|_{\K}^2 \,ds\Big)
\end{align}
Regarding the duality argument used to estimate $\E\int_{T-\delta}^T|Q(s)|^2_{\mathcal{K}} ds $ the only thing to check is that \eqref{stima5}  still holds:
\begin{align}\label{stima5Lyap}
&\Big|  \E \int_{T-\delta} ^T \langle X^n(s), J_nS(s)J_n\rangle_{L_2(H)}  \, ds \Big| =
\Big|  \E \int_{T-\delta} ^T \sum_{k\geq 1} \langle X^n(s)e_k, J_nS(s)J_ne_k\rangle_{H}  \, ds \Big|
\\ & \leq  \left(\int_{T-\delta} ^T \E  \sum_{k\geq 1} \lambda_k^{2\rho} | X^n(s) e_k|^2  \,ds\right)^{1/2} 
  \left(\int_{T-\delta} ^T \E  \sum_{k\geq 1} \lambda_k^{-2\rho} | J_nS(s)J_ne_k|^2  \,ds\right)^{1/2}  \nonumber \\ &\leq \left(\delta \int_{T-\delta} ^T
|Q^n(s)|^2_{\K} \,ds\right)^{1/2} \left(  \int_{T-\delta} ^T \E \sum_{k\geq 1} \lambda_k^{-2\rho} | S(s)e_k|^2 \,ds\right)^{1/2}\nonumber \\ &\leq
 \delta ^{1/2}\left(\int_{T-\delta} ^T
|Q^n(s)|^2_{\K} \,ds\right)^{1/2}\left( \int_{T-\delta} ^T\nonumber
|S(s)|^2_{\K} \,ds\right)^{1/2}
\\\nonumber& \leq  2\delta \E \int_{T-\delta}^T  |S(s)|^2_{L(H)} ds+ \frac{1}{8}\E \int_{T-\delta}^T  |Q^n(s)|^2_\K  \,ds
\end{align}
Thus we deduce again \eqref{stimaQnunif}, that together with \eqref{StimaPUnifLyap}  leads to prove \eqref{aprioriLyap}.

We also have that
\begin{theorem}\label{lyapunovMax.teo}
Under  assumptions $\rm{A1)-A2)-A3')}$ equation
\eqref{BSDEoperatorvalued} has a unique mild solution $(P,Q)$.
\end{theorem}
\Dim
The only thing to check is that following  still hold:
\begin{equation} \label{convSLyap} \lim_{n,m \to +\infty}\E \int_{T-\delta}^T |J_nS(s)J_n -J_mS(s)J_m|^2_\K\,ds =0
\end{equation}
Recalling that $ e_k \in V$,  for every $ k \geq 1$, we have:
\begin{equation} \lim_{n,m \to +\infty}|J_nS(s)J_n e_k -J_mS(s)J_m e_k|^2_\K\,ds =0, \qquad \forall k \geq 1.
\end{equation}
Moreover
\begin{equation}\E \int_{T-\delta}^T |J_nS(s)J_n -J_mS(s)J_m|^2_\K\,ds \leq 2 \E \int_{T-\delta}^T |S(s)|^2_\K\,ds
\end{equation}
Thus by Dominated Convergence Theorem we get that \eqref{convSLyap}.
The rest of the proof follows then identically as Theorem \ref{lyapunov.teo}. \finedim

\begin{example}
Notice that in the mentioned papers \cite{Fu-Hu-Te-1} and \cite{Fu-Hu-Te-2}, $H=L^2([0,1])$ and the operator $S(t)$ is the multiplication operator by an adapted  stochastic random field $H:  (\Omega\times[0,T]\times[0,1]) \rightarrow \mathbb{R}$ namely
$$ [S(t) e](\xi)=H(t,x)e(x),\;  \forall e\in L^{\infty}([0,1]),\; \forall x\in [0,1], \quad \hbox{ with }  \quad \E\int_0^T\int _0^1 H(t,x)^2 dt dx < +\infty$$
(notice that in thes case $S(t)$ is not even defined on the whole $H$).

Moreover the infinitesimal generator $A$ is the realization of the Laplacian in $L^2([0,1])$  with Dirichlet boundary conditions.

Thus we have, choosing the basis $\{e_m\}_{m\in \mathbb{N}}$, of eigenvectors of $A$:
\begin{itemize}
\item[(a)] $\sup_{m \geq 1 } |e_m|_{L^\infty([0,1])} < \infty$.
\item [(b)] $S$ is self-adjoint and $|S(s)e_m|^2= \int_0^1 (H ^2(s,x)e_m^2(x) \,dx \leq \sup_{m \geq 1 } |e_m|^2_{L^\infty([0,1])} |H(s,\cdot)|_{L^2([0,1])}$
\item[(c)] $|S(s)|_{\K}= \sum_{k \geq 1} \sum_{m \geq 1}  \lambda_m^{-2\rho}|\langle S(s) e_k, e_m \rangle _{L^2([0,1])}|^2=  \sum_{m \geq 1}  \lambda_m^{-2\rho}  \sum_{k \geq 1}|\langle e_k, S (s)e_m \rangle _{L^2([0,1])}|^2  $  $ =  \sum_{m \geq 1}  \lambda_m^{-2\rho} |S(s)e_m|^2_{L^2([0,1])}\leq  |H^2(s,\cdot)|_{L^2([0,1])} \sum_{m \geq 1}  \lambda_m^{-2\rho}\leq \hbox{cost} |H^2(s,\cdot)|_{L^2([0,1])} $
\end{itemize}
and Assumptions ${\rm A1)\, A2)\, A3') }$ hold.
\end{example}

\end{document}